\documentclass[a4paper,10pt]{amsart} 
\usepackage{mathptmx}
\usepackage{amsfonts,amssymb} 
\usepackage{amsmath} 
\usepackage{multicol}

\renewcommand{\setminus}{\smallsetminus}

\newcommand{\iso}{\cong}

\newcommand{\Z}{{\mathbb Z}}

\newcommand{\fof}{\textit{ff}}

\newtheorem{theorem}{Theorem}[section]

\newtheorem{corollary}[theorem]{Corollary}

\newtheorem{definition}[theorem]{Definition}
\newtheorem{example}[theorem]{Example}

\newtheorem{lemma}[theorem]{Lemma}

\newtheorem{proposition}[theorem]{Proposition}

\newtheorem{statement}[theorem]{}

\newcommand{\F}{\mathbb F} \newcommand{\rad}{\operatorname{Rad}}
\renewcommand{\ker}{\operatorname{Ker}} 

\title{Invariant rings of orthogonal groups over $\F_2$} 
\author{P. H. Kropholler}
\author{S. Mohseni Rajaei}
\author{J. Segal} 
\address{Dept of Mathematics, University of Glasgow, University Gardens, Glasgow G12 8QW}
\email{p.h.kropholler@maths.gla.ac.uk} 
\address{Departement of Mathematics, Azzahra University, Vanak, Tehran,Iran}
\email{rajaei@azzahra.ac.ir} 
\address{Stegem\"uhlenweg 70, 37083 G\"ottingen}
\email{joel@berlin.com} 

\begin{document}

\begin{abstract}
We determine the rings of invariants $S^G$ where $S$ is the symmetric algebra on the dual of a vector space $V$ over $\F_2$ and $G$ is the orthogonal group preserving a non-singular quadratic form on $V$. The invariant ring is shown to have a presentation in which the difference between the number of generators and the number of relations is equal to the minimum possibility, namely $\dim V$, and it is shown to be a complete intersection. In particular, the rings of invariants computed here are all Gorenstein and hence Cohen-Macaulay.
\end{abstract}

\maketitle 

\tableofcontents

\section*{Acknowledgements}

The first author wishes to acknowledge his collaboration with David Carlisle during the mid nineteen-eighties when the calculations described here were first studied.

\section{Strategy for calculating invariants}

Let $S$ be the symmetric algebra on $V^*$, the dual of a finite dimensional vector space $V$. For any finite subgroup $G$ of $GL(V)$ we can consider the invariant ring $S^G$. This article concerns explicit calculations of $S^G$ which are described in \S6. In very broad outline, the strategy for doing calculations comprises the following steps:
\begin{itemize}
\item Find some reasonably large  but finite collection of invariants of $G$ using a variety of methods.
\item Consider the subring of $T$ which they generate. If we found enough invariants in the first step then $T$ will be $S^G$ and we are done. The strategy cannot be doomed to failure because $S^G$ is a finitely generated ring.
\item Prove that $T=S^G$.
\end{itemize}
The last step may fail. If it does, then we hope to discover new invariants. We throw these in to the generating set, enlarge $T$ and try again. 
The proof is easy enough:
\begin{itemize}
\item Show that $S$ is integral over $T$. 
\item Show that $T$ has the right field of fractions.
\item Show that $T$ is integrally closed.
\end{itemize}
Only the third item here causes any real concern. In fact, we shall be working on examples which are known {\em a priori} to be unique factorization domains and we'll establish this for our $T$ as the route to integral closure. 
Over $\F_2$ invariant rings are always unique factorization domains by a result of Nakajima, see Corollary 3.9.3 of \cite{benson}. This is not true in odd characteristic as illustrated by the very simple example 1 of Chapter 1 of \cite{smith}. It is also not true in general for fields of characteristic $2$ which have more than two elements.
The following elementary result turns out to be decisive in every case considered in this paper:
\begin{proposition}
Let $R$ be a commutative ring and suppose that $\alpha,\beta$ is a regular sequence in $R$ such that
\begin{enumerate}
\item the localizations $R[\alpha^{-1}]$ is a unique factorization domains;
\item $\alpha$ generates a prime ideal in the ring $R[\beta^{-1}]$.
\end{enumerate}
Then $R$ is a unique factorization domain.
\end{proposition}
\begin{proof}
First, $\alpha$ being a non-zero divisor, the map $R\to R[\alpha^{-1}]$ is injective and so we know that $R$ is a domain. Since  $\alpha$ is prime in $R[\beta^{-1}]$, the quotient $R[\beta^{-1}]/\alpha R[\beta^{-1}]$ is a domain. Now $\beta$ is a non-zero-divisor modulo $\alpha$ and so the map $R/\alpha R\to R[\beta^{-1}]/\alpha R[\beta^{-1}]$ is injective. Therefore $\alpha$ generates a prime ideal in $R$, and this together with the fact that $R[\alpha^{-1}]$ is a unique factorization domain implies that $R$ is a unique factorization domain.
\end{proof}

\section{Introduction} 

In this paper we shall focus attention on an odd-dimensional vector
space $V$ over the field $\F_2$ of two elements which is endowed with a
non-singular quadratic form $\xi_0$. We write $S=S(V^*)$ for the symmetric
algebra on the dual $V^*$ of $V$. The symmetric algebra is the
polynomial ring in any chosen basis of $V^*$, and it inherits a natural
action of $GL(V)$. The orthogonal group of automorphisms of $V$ which
preserve $\xi_0$ is denoted by $O(V)$. Our main objective is to compute
the ring of invariants of $O(V)$. The results extend work \cite{rajaei,rajaei2} of the second author on the rational invariants of orthogonal groups.

The reader familiar with quadratic forms in characteristic $2$, the associated finite orthogonal groups, the rudiments of the mod $2$ Steenrod algebra and the Dickson invariants for the general linear group may now wish to skip straight to \S6 for a statement of results.

We begin with some remarks which apply to any non-zero vector space $V$ over any
field $K$. A quadratic form $q$ is a function $$q:V\to K$$ which
satisfies the two conditions \begin{itemize} \item the polarization
$$b:V\times V\to K$$ defined by $$b(u,v)=q(u+v)-q(u)-q(v)$$ is
bilinear; and \item for all scalars $\lambda$ and all $v\in V$,
$$q(\lambda v)=\lambda^2q(v).$$ \end{itemize} 
The polarization $b$ is always a symmetric form.
If the characteristic of $K$ is not $2$ then the quadratic form $q$ can be recovered from its polarization by the formula $q(v)=\frac12b(v,v)$ and there is a bijective correspondence between quadratic forms and symmetric bilinear forms. If $K$ has characteristic $2$ then the polarization is an alternating form from which the quadratic form cannot be recovered. If $K=\F_2$, the case of interest in this paper, then the definition of quadratic form simplifies: a function $q:V\to\F_2$ such that
\begin{itemize} \item the polarization
$$b:V\times V\to \F_2$$ defined by $$b(u,v)=q(u+v)+q(u)+q(v)$$ is
bilinear. \end{itemize}  In this case, every alternating form arises as the polarization of $2^m$ different quadratic forms, where $m=\dim V$, and the symmetric forms which are not alternating never arise as polarizations. Thus polarization yields a map between quadratic forms and symmetric bilinear forms which is neither injective nor surjective.

%%%%%%%%%%%%%%%%%%%%%%%%%%%%%%%%%%%%%%%%%%%%%%%%%%%%%%%%%%%%%%%%%%%

\section{The Steenrod Algebra and Chern polynomials} 

Henceforth we assume that $V$ is a vector space over $\F_2$. The symmetric
algebra $S=S(V^*)$ is naturally isomorphic to the cohomology ring
$H^*(BV,\F_2)$ of the classifying space $BV$ of the additive group $V$, drawing attention to the fact that $S$ admits an unstable
action of the Steenrod Algebra ${\mathcal A}_2$. The reader is referred
to the books \cite{smith,steenrod-epstein} for details about the
Steenrod Algebra. 

Here the matter is simple enough. The Steenrod algebra is an $\F_2$-algebra generated 
by elements $Sq^i$, for $i\ge0$, called {\em Steenrod squares}. 
$Sq^i$ is homogeneous of degree $i$.
\begin{statement}
The action on $S$ is determined by the following facts: 
\begin{itemize} 
\item $Sq^0$ acts as the identity operation on $S$; 
\item $Sq^1$ acts as a derivation on $S$; 
\item for all $x\in V^*$,
$Sq^1(x)=x^2$ and $Sq^n(x)=0$ for $n\ge2$; 
\item the Cartan formula holds: for all $s$ and $t$ in $S$, $$Sq^n(st)=\sum_{i+j=n}\left(Sq^is\right)\left(Sq^jt\right).$$
\item for any
homogeneous element $s$ of $S$ of degree $d$, $Sq^d(s)=s^2$ and
$Sq^j(s)=0$ if $j>d$. 
\item The {\em total} Steenrod operation $Sq^\bullet:=Sq^0+Sq^1+Sq^2+\cdots$ acts as a ring endomorphism of $S$.
\end{itemize} 
\end{statement}

Now suppose that $\mathfrak{S}$ is a non-empty subset of $V^*$ which contains $d$
elements. The {\em Chern polynomial} associated to $\mathfrak{S}$ is the
polynomial $$\prod_{x\in\mathfrak{S}}(X+x).$$ Let's write $f_i$ for the
coefficient of $X^{d-i}$ so that
$$\prod_{x\in\mathfrak{S}}(X+x)=f_0X^d+f_1X^{d-1}+\dots+f_d.$$ Then it is easy to
see that 
\begin{lemma}\label{insight}
For each $i$ in the range $0\le i\le d$, $$Sq^i(f_d)=f_df_i.$$ 
\end{lemma}
This is a special case of the Wu formulae \cite{koch,stong,wu} for the action of the Steenrod algebra on the cohomology ring $H^*(BO,\F_2)$ of the classifying space for real vector bundles which carries the generic Stiefel-Whitney classes. Arguably the name ``Stiefel-Whitney  polynomial'' would be more appropriate in this paper than our choice: Chern polynomial. On the other hand, in modular invariant theory the similarity between the characteristic $2$ theory and the odd characteristic theory is close and we stick with the name Chern polynomial.

%%%%%%%%%%%%%%%%%%%%%%%%%%%%%%%%%%%%%%%%%%%%%%%%%%%%%%%%%%%%%%%%%%%%%%%

\section{Quadratic forms over $\F_2$}\label{quadforms}

Each element of the symmetric algebra $S$ determines a function from $V$ to $\F_2$ and in this
way the homogeneous elements of degree two in $S$ determine quadratic
forms on $V$. Conveniently, it is the case that this correspondence
between $S_2$ and the set of quadratic forms on $V$ is a bijection. 
We identify quadratic forms with the corresponding elements of $S_2$.

\begin{lemma} Let $q$ and $q'$ be quadratic forms on $V$. Then the
following are equivalent: \begin{enumerate} \item $q$ and $q'$ have the
same polarization; \item $Sq^1(q)=Sq^1(q')$; \item $q+q'=x^2$ for some
$x\in V^*$. \end{enumerate} 
\end{lemma} 
\begin{proof} The details of
this easy lemma are left to the reader. Note that the Steenrod operation
$Sq^1$ is determined by virtue of being a derivation such that
$Sq^1(x)=x^2$ for all $x\in V^*$. 
\end{proof} 

Two alternating forms $b$ and $b'$ on $V$ are {\em equivalent} iff there exists $g\in GL(V)$ such that
$b'(u,v)=b(gu,gv)$ for all $u$, $v$.
Alternating forms are determined up to equivalence by their rank, and
their rank is always even. Let $b$ be an alternating form on $V$. The
radical $\rad(b)$ of $b$ is defined to be $$\{v\in V;\ b(v,\ \ )=0\}.$$
If $q$ is a quadratic form which polarizes to $b$ then the radical
$\rad(q)$ of $q$ is defined to be $$\{v\in\rad(b);\ q(v)=0\}.$$ Since
the restriction of $q$ to $\rad(b)$ is a linear functional, one finds
that $\rad(q)$ is either equal to $\rad(b)$ or has codimension $1$ in
$\rad(b)$. A quadratic form $q$ is called {\em non-singular} if and only if $\rad(q)=0$. 

We consider quadratic forms always in the presence of a fixed alternating form to which they polarize. 
We shall use the term {\em symplectic space} to refer to a finite dimensional vector space endowed with an alternating form of maximum possible rank. The group of automorphisms of a symplectic space is called the {\em symplectic group}. Non-singular quadratic forms live on symplectic spaces.
On a symplectic space, we say that two quadratic forms $q$ and $q'$ are {\em equivalent} iff there exists $g$ in the symplectic group such that $q'(v)=q(gv)$ for all $v$. 

On a non-zero even dimensional symplectic space there are two types of non-singular
quadratic form up to equivalence, called $+$type and $-$type. In
dimension $2n\ge2$ 
\begin{itemize} \item $2^{2n-1}+2^{n-1}$ of these
forms have $+$type, 
\item $2^{2n-1}-2^{n-1}$ of these forms have $-$type. 
\end{itemize} These quadratic forms are also classified by the Arf invariant which is determined by Browder's democracy: the Arf invariant is the value of the quadratic form taken by a majority of vectors. The forms of $+$type have Arf invariant $0$ and the forms of $-$type have Arf invariant $1$.

If $V$ is an odd dimensional symplectic space, then there is only one kind of
non-singular quadratic form up to equivalence. If $\xi_0$ is such a form
and $b$ is its polarization, then each form having the same polarization
is equal to $\xi_0+x^2$ for some $x\in V^*$ and there are three kinds:
the non-singular forms (all equivalent to $\xi_0$), the
singular forms of $+$type, and the singular forms of $-$type. 

In this case, the polarization $b$ is degenerate and its radical contains a vector $e_0\ne0$. Correspondingly there is a subspace $U^*$ of $V^*$ of codimension one and the forms of non-singular type are exactly the forms $\xi_0+x^2$ for $x\in U^*$.

If $\dim
V=2n+1\ge3$ then 
\begin{itemize} 
\item $2^{2n}$ of these forms are
non-singular, and each is equal to $\xi_0+x^2$ for some $x\in V^*$ such
that $\ker x\supseteq\rad(b)$; 
\item $2^{2n-1}+2^{n-1}$ of these forms
have $+$type, and each is equal to $\xi_0+x^2$ for certain $x\in V^*$
such that $\ker x\cap\rad(b)=0$; 
\item $2^{2n-1}-2^{n-1}$ of these forms
have $-$type, and each is equal to $\xi_0+x^2$ for certain $x\in V^*$
such that $\ker x\cap\rad(b)=0$. 
\end{itemize} 

Note that for a $+$ or $-$type form $q$, this implies that 
$q=\xi_0+x^2$ 
is actually in $S(U^*)$. Note also that the Arf invariant is not defined for the non-singular forms.

%%%%%%%%%%%%%%%%%%%%%%%%%%%%%%%%%%%%%%%%%%%%%%%%%%%%%%%%%%%%%%%%%%%%%%%

\section{Definitions} \label{section definitions}

In this section we state our definitions of symplectic and orthogonal
groups and those of their representations that we consider.

Let $n$ be a positive integer. Henceforth we suppose that $V$ has dimension $2n+1$ over $\F_2$ and that $\xi_0$ is a
non-singular quadratic form on $V$. Let $b$ denote the polarization of
$\xi_0$. The radical of $b$ is one dimensional. Choose a basis
$e_0,\dots,e_{2n}$ of $V$ where $e_0$ is the non-zero vector in the
$\rad(b)$. 

The reader is referred to Cameron's notes \cite{cameron} for the
background to the following lemma and definition. 
\begin{lemma} \label{standard basis} It is possible to choose the $e_i$
for $i\ge1$ so that the matrix $B$ with $(i,j)$-entry
$B_{i,j}=b(e_i,e_j)$ is 
$$\left(\begin{matrix} 0&&&&&&&\\ &0&1&&&&&\\ &1&0&&&&&\\ &&&0&1&&&\\
&&&1&0&&&\\ &&&&&\ddots&&\\ &&&&&&0&1\\ &&&&&&1&0\\
\end{matrix}\right)$$ filled out with zeroes. \end{lemma} 
We write $B_0$ for the non-singular alternating matrix obtained by
omitting the first row and column of $B$. 
\begin{definition} \label{groups}  
\begin{enumerate} 
\item $O(V)$
denotes the (orthogonal) group of automorphisms of $V$ which preserve
the quadratic form $\xi_0$. \item $Sp(V)$ denotes the (symplectic) group
of automorphisms of $V$ which preserve the alternating form $b$.
$$O(V)\subset Sp(V).$$ \item $U$ denotes the quotient $V/\langle
e_0\rangle$. This space inherits the alternating form, but it does not
inherit any natural quadratic form. \item $Sp(U)$ denotes the
(symplectic) group of automorphisms of $U$ which preserve the inherited
alternating form. 
\end{enumerate} 
\end{definition} 

Let $x_0,\dots,x_{2n}$ be the basis of $V^*$ which is dual to our chosen
basis of $V$. We remark that when the basis $e_i$ is chosen in
accordance with Lemma \ref{standard basis}, then the quadratic form is
given by $$\xi_0=x_0^2+x_1x_2+x_3x_4+\dots+x_{2n-1}x_{2n}.$$ The
canonical surjection $V\to U$ induces an injection $U^*\to V^*$. We
identify $U^*$ with its image in $V^*$: thus $U^*$ is the subspace of
$V^*$ spanned by $x_1,\dots,x_{2n}$. The symmetric algebra on $U^*$ is
the subring of $S$ generated by $x_1,\dots,x_{2n}$. Every symplectic
automorphism of $V$ induces a symplectic automorphism of $U$, and every
symplectic automorphism of $U$ arises this way. In this way there is a
surjective homomorphism $$Sp(V)\to Sp(U).$$ The kernel of this
homomorphism consists of transvections: it is the elementary abelian
$2$-group of rank $2n$ comprising the linear automorphisms of $V$ which
fix $e_0$ and induce trivial action on $U$, and can be naturally
identified with $\hom(U,\langle e_0\rangle)\iso U^*$. The homomorphism
between the symplectic groups restricts to an isomorphism $$O(V)\iso
Sp(U).$$ 

\begin{definition} The sequence $\xi_1,\xi_2,\xi_3,\dots$ is defined
recursively by $$\xi_{n}=Sq^{2^{n-1}}(\xi_{n-1}).$$ \end{definition} 

When the basis $e_i$ is chosen in accordance with Lemma \ref{standard
basis}, then 
$$\xi_j=x_1^{2^j}x_2+x_1x_2^{2^j}+x_3^{2^j}x_4+x_3x_4^{2^j}+\dots+x_{2n-
1}^{2^j}x_{2n}+x_{2n-1}x_{2n}^{2^j}$$ for each $j\ge1$.

In general, for $i\ge1$, each $\xi_i$ belongs to the symmetric algebra
on $U^*$ (i.e. it does not involve $x_0$) and is an invariant of
$Sp(U)$. For each $i\ge0$, $\xi_i$ has degree $2^i+1$. The following results are important: 

\begin{lemma}\label{squarexi}%
%\begin{multicols}{2}%
%\begin{eqnarray*}Sq^0(\xi_0)&=&\xi_0\\ Sq^1(\xi_0)&=&\xi_1\\
%Sq^2(\xi_0)&=&\xi_0^2\\ Sq^j(\xi_0)&=&0\text{\ \ in all other cases.}
%\end{eqnarray*}
%\begin{eqnarray*}
%&&\\
%Sq^0(\xi_1)&=&\xi_1\\ Sq^2(\xi_1)&=&\xi_2\\ Sq^3(\xi_1)&=&\xi_1^2\\
%Sq^j(\xi_1)&=&0\text{\ \ in all other cases.}
%\end{eqnarray*} 
%\end{multicols}
%and for
%all $i\ge2$, 
%\begin{eqnarray*} Sq^0(\xi_i)&=&\xi_i\\
%Sq^1(\xi_i)&=&\xi_{i-1}^2\\ Sq^{2^i}(\xi_i)&=&\xi_{i+1}\\
%Sq^{2^i+1}(\xi_i)&=&\xi_i^2\\ Sq^j(\xi_i)&=&0\text{\ \ in all other
%cases.} 
%\end{eqnarray*}\ 
%\end{lemma} 
%These results can be conveniently summarized using the total operation:
\begin{eqnarray*}
Sq^\bullet(\xi_0)&=&\xi_0+\xi_1+\xi_0^2,\\
Sq^\bullet(\xi_1)&=&\xi_1+\xi_2+\xi_1^2,\\
Sq^\bullet(\xi_i)&=&\xi_i+\xi_{i-1}^2+\xi_{i+1}+\xi_i^2\qquad(i\ge2).
\end{eqnarray*}
\end{lemma}
\begin{proof}
It is easy given that $Sq^\bullet$ is a ring homomorphism and $Sq^\bullet x=x+x^2$ for $x\in V^*$.
\end{proof}
\begin{corollary} For any $j$ and
any $m\ge0$, \begin{itemize} \item
$Sq^j\left(\F_2[\xi_0,\dots,\xi_m]\right)\subseteq\F_2[\xi_0,\dots,
\xi_{m+1}]$; \item
$Sq^j\left(\F_2[\xi_1,\dots,\xi_m]\right)\subseteq\F_2[\xi_1,\dots,
\xi_{m+1}]$. \end{itemize} \end{corollary} 

\begin{definition}\label{dickthm}
Elements $c_0,\dots, c_{2n}$ of $S(U^*)$ are defined to
be the unique elements of $S$ such that $$\prod_{x\in
U^*}(X+x)=\sum_{j=0}^m c_{j}X^{2^j}.$$ These are the Dickson invariants:
$$S(U^*)^{GL(U)}=\F_2[c_{2n-1},\dots,c_0].$$ We write $D(X)$ for the
Dickson polynomial. 
\end{definition} 

Notice that the Dickson polynomial is the Chern polynomial associated to
the subset $\mathfrak{S}:=U^*$. Using Lemma \ref{insight} we have that

\begin{lemma} \label{steenroddickson} 
For $0\le i\le2n$,
$$Sq^{2^{2n}-2^i}(c_0)=c_0c_i$$ and $Sq^j(c_0)$ is zero in all other
cases. 
\end{lemma} 

Dickson's original paper
\cite{dickson} introduced these invariants: for a more modern treatment see Wilkerson \cite{wilkerson}. Crucially, $c_0D(X)$ is equal to the
determinant $$\left| \begin{matrix} X&x_1&x_2&x_3&\dots&x_{2n}\\\\
X^2&x_1^2&x_2^2&x_3^2&\dots&x_{2n}^2\\\\
X^4&x_1^4&x_2^4&x_3^4&\dots&x_{2n}^4\\\\
X^8&x_1^8&x_2^8&x_3^8&\dots&x_{2n}^8\\\\
\vdots&\vdots&\vdots&\vdots&\ddots&\vdots\\\\
X^{2^{2n}}&x_1^{2^{2n}}&x_2^{2^{2n}}&x_3^{2^{2n}}&\dots&x_{2n}^{2^{2n}}\
\ \end{matrix} \right|.$$ 
Let $C_0$ denote the matrix 
$$\left(
\begin{matrix} x_1&x_2&x_3&\dots&x_{2n}\\\\
x_1^2&x_2^2&x_3^2&\dots&x_{2n}^2\\\\
x_1^4&x_2^4&x_3^4&\dots&x_{2n}^4\\\\
x_1^8&x_2^8&x_3^8&\dots&x_{2n}^8\\\\
\vdots&\vdots&\vdots&\ddots&\vdots\\\\
x_1^{2^{2n-1}}&x_2^{2^{2n-1}}&x_3^{2^{2n-1}}&\dots&x_{2n}^{2^{2n-1}}\\
\end{matrix} \right).$$ 
Then $C_0$ has determinant $c_0$ and the matrix
equation below simply expresses the fact that $D(x_i)$ vanishes for each
$i$. 

\begin{lemma} \label{prerelations} 
$$C_0^T\left( \begin{matrix}
c_0\\ c_1\\ c_2\\ \vdots\\ c_{2n-1} \end{matrix} \right)=\left(
\begin{matrix} x_1^{2^{2n}}\\ x_2^{2^{2n}}\\ x_3^{2^{2n}}\\ \vdots\\
x_{2n}^{2^{2n}} \end{matrix} \right).$$ \end{lemma} For later use, we
write $C$ for the $(2n+1)\times2n$-matrix $$\left( \begin{matrix}
x_1&x_2&x_3&\dots&x_{2n}\\\\ x_1^2&x_2^2&x_3^2&\dots&x_{2n}^2\\\\
x_1^4&x_2^4&x_3^4&\dots&x_{2n}^4\\\\
x_1^8&x_2^8&x_3^8&\dots&x_{2n}^8\\\\
\vdots&\vdots&\vdots&\ddots&\vdots\\\\
x_1^{2^{2n-1}}&x_2^{2^{2n-1}}&x_3^{2^{2n-1}}&\dots&x_{2n}^{2^{2n-1}}\\\\
x_1^{2^{2n}}&x_2^{2^{2n}}&x_3^{2^{2n}}&\dots&x_{2n}^{2^{2n}}\\
\end{matrix} \right),$$ and we write $\widehat C$ for the
$(2n+1)\times(2n+1)$-matrix $$\left( \begin{matrix}
x_0&x_1&x_2&x_3&\dots&x_{2n}\\\\
x_0^2&x_1^2&x_2^2&x_3^2&\dots&x_{2n}^2\\\\
x_0^4&x_1^4&x_2^4&x_3^4&\dots&x_{2n}^4\\\\
x_0^8&x_1^8&x_2^8&x_3^8&\dots&x_{2n}^8\\\\
\vdots&\vdots&\vdots&\vdots&\ddots&\vdots\\\\
x_0^{2^{2n}}&x_1^{2^{2n}}&x_2^{2^{2n}}&x_3^{2^{2n}}&\dots&x_{2n}^{2^{2n}
}\\ \end{matrix} \right).$$ 

\begin{lemma} \label{dick} 
Suppose that $f_0,\dots,f_{2n}$ are elements
of $S$ with the property that for all $i$ in the range $1$ to $2n$,
$$f_{2n}x_i^{2^{2n}}+f_{2n-1}x_i^{2^{2n-
1}}+\dots+f_{2}x_i^{4}+f_{1}x_i^{2}+f_0x_i=0.$$ Then $$f_j=f_{2n}c_j$$
for all $j$. \end{lemma} 

\begin{proof} The polynomial
$$f(X)=f_{2n}X^{2^{2n}}+f_{2n-1}X^{2^{2n-
1}}+\dots+f_{2}X^{4}+f_{1}X^{2}+f_0X$$ vanishes on all $x_i$. The
additivity of the Frobenius map enables us to conclude that $f(X)$
vanishes on any linear combination of the $x_i$: that is, $f(X)$
vanishes on $U^*$. Therefore $f(X)$ is divisible by the Dickson
polynomial, and the result follows. \end{proof} 

%%%%%%%%%%%%%%%%%%%%%%%%%%%%%%%%%%%%%%%%%%%%%%%%%%%

\section{Statement of Results for the Orthogonal Groups}

We give a summary of the conclusions of the calculations.

Assume that $n\ge2$ and that $V$ is a $(2n+1)$-dimensional $\F_2$-space endowed with a non-singular quadratic form $\xi_0$. The cases $n\le 1$ will be treated later and are quite elementary by comparison.
\begin{theorem}
Let $T^\dagger$ be an abstract polynomial ring on generators $$\xi_0,\dots,\xi_{2n-2},d_{2n-1},\dots,d_n$$ where $\xi_i$ has degree $2^i+1$ and $d_j$ has degree $2^{2n-1}-2^{j-1}$. Define the degree preserving map
$$T^\dagger\to S$$ by sending $\xi_0$ to the quadratic form of the same name, sending $\xi_i$ to $Sq^{2^{i-1}}Sq^{2^{i-2}}\dots Sq^1\xi_0$ and sending $d_j$ to the symmetric polynomial of the same degree in the set of vectors in $V^*$ of $-$type (i.e. the elements $x$ of $V^*$ such that $\xi_0+x^2$ has $-$type in $S$). Then the image of $T^\dagger$ in $S$ is equal to the ring $S^{O(V)}$ of invariants of the orthogonal group of automorphisms preserving $\xi_0$. The kernel of the map is generated by a regular sequence of $n-2$ elements which are homogeneous of degrees $2^{2n-1}+2^j$ for $1\le j\le n-2$. 

The $n-2$ relations can be expressed in matrix form as the sum of four column vectors not all of which are easily described at this stage. When they are expressed this way, we use an additional (redundant) generator $\xi_{2n-1}$ and impose an additional relation at the beginning which amounts to an expression for $\xi_{2n-1}$ in terms of the chosen generators:
$$\mathbf{S}_{n-1}
\left(\begin{matrix}d_n\\\vdots\\d_{2n-2}\end{matrix}\right)+
\left(\begin{matrix}\xi_{2n-1}\\\vdots\\\xi_{n+1}^{2^{n-2}}\end{matrix}\right)+
\sqrt{\left(\mathbf{L}'_n\mathbf{E}_n+\left(\mathbf{L}'_n\mathbf{K}_n+
\mathbf{R}'_n\right)\mathbf{F}_n\right)}+\mathbf{G}_{n-1}d_{2n-1}.$$
Here, the matrices and vectors 
$\mathbf{S}_{n-1},
\mathbf{L}'_n,\mathbf{R}'_n,\mathbf{E}_n,\mathbf{F}_n,\mathbf{G}_{n-1}$ are defined subsequently and all involve only polynomials in the $\xi$'s. 
Matrices $\mathbf{L}_n,\mathbf{R}_n$ are defined in the next section and the symbol ${}'$ above indicates the matrices obtained from these by deleting the first row. The matrix $\mathbf{K}_n$ and column vector $\mathbf{E}_n$ are also introduced in the next section. The column vector $\mathbf{F}_n$ arises as part of a connection (Lemma \ref{frcd}) between the symmetric polynomial invariants determined by $d_{2n-1},\dots,d_n$ and the Dickson invariants $c_{2n-1},\dots,c_n$ for the $2n$-dimensional quotient of $V$ which inherits a natural alternating form. The matrix $\mathbf{S}_{n-1}$ has determinant equal to $\Lambda_{2n-2}$, a certain polynomial which is intimately related to the Dickson algebra for a vector space of dimension $2n-2$. The $\sqrt{}$ symbol here is used to indicate the matrix obtained by replaced each entry of the matrix to which it is applied by its square root. We shall see that every entry of $\left(\mathbf{L}'_n\mathbf{E}_n+\left(\mathbf{L}'_n\mathbf{K}_n+
\mathbf{R}'_n\right)\mathbf{F}_n\right)$ is a square in $T^\dagger$.

The second of the four column vectors has $i${\em th} entry $\xi_{2n-i}^{2^{i-1}}$ and the $i${\em th} relation can be interpreted as saying that $\xi_{2n-i}^{2^{i-1}}$ can be expressed in terms of other generators and lower powers of $\xi_{2n-i}$. This relation cannot be deduced from any of the other relations. 

The regular sequence of relations in $T^\dagger$ can be extended to a regular sequence of length equal to $3n-1$, the Krull dimension of $T^\dagger$, by taking the further $2n+1$ elements
$$\xi_0,\dots,\xi_n,d_{2n-1},\dots,d_n.$$
The ring of invariants $S^{O(V)}$ is a complete intersection as described in \S21 of \cite{mats}.
\end{theorem}

Notice that we use symmetric polynomials arising from vectors of $-$type rather than $+$type to describe the invariants $d_j$. The reason for this is purely pragmatic: there are fewer vectors of $-$type.

In general for any symmetric algebra $S$ on the dual of a finite vector space we know the following.
The Hilbert series of any ring of invariants $S^G$ has a Laurent expansion about $t=1$ which begins with 
$$\frac{1}{|G|}\frac{1}{(1-t)^m}+\frac{r}{2|G|}\frac{1}{(1-t)^{m-1}}+\dots$$
where $m$ is the Krull dimension of $S$ (i.e. the dimension of $V$) and $r$ is determined by the reflections (i.e. elements fixing a hyperplane in $V$ pointwise) using a ramification formula. 
The interpretation of the second coefficient in terms of reflections is the content of the Benson--Crawley-Boevey--Neeman theorem. This was first proved by Benson and Crawley-Boevey, see \cite{BCB}. Subsequently, Neeman published a line of reasoning which uses the Riemann--Roch theorem, \cite{Neeman}. Evidence that the Riemann--Roch theorem is involved was observed much earlier in unpublished work \cite{voloch} of Felipe Voloch. 

Since our invariant ring is presented by means of a regular sequence the Hilbert series is very simply determined. The Hilbert series of the ring $T^\dagger$ is
$$\frac{1}{(1-t^2)(1-t^3)\dots(1-t^{2^{2n-2}+1})
\cdot(1-t^{2^{2n-1}-2^{2n-2}})\dots(1-t^{2^{2n-1}-2^{n-1}})}.$$ 
There is a contribution $(1-t^a)$ in the denominator for a generator of degree $a$. Each time a relation of degree $b$ is imposed we simply multiply this Hilbert series by $(1-t^b)$ because the relation is a non-zero-divisor modulo its predecessors.
We can therefore draw the following conclusions:
\begin{corollary}
The ring of invariants $S^{O(V)}$ has Hilbert polynomial
$$\frac{(1-t^{2^{2n-1}+2})\dots(1-t^{2^{2n-1}+2^{n-2}})}
{(1-t^2)(1-t^3)\dots(1-t^{2^{2n-2}+1})
\cdot(1-t^{2^{2n-1}-2^{2n-2}})\dots(1-t^{2^{2n-1}-2^{n-1}})}.$$
\end{corollary}
\begin{proof}
There is a contribution $(1-t^b)$ in the numerator for a relation of degree $b$.
\end{proof}
Noting that the Laurent power series expansion of 
$$\frac{\prod_{i=1}^k(1-t^{b_i})}{\prod_{j=1}^\ell(1-t^{a_i})}$$ about $t=1$ begins
$$\frac{\prod b_j}{\prod a_i}\left(\frac{1}{(1-t)^m}+
\frac{\sum(a_i-1)-\sum(b_j-1)}{2}\frac{1}{(1-t)^{m-1}}+\dots\right)$$ where $m=\ell-k$, we can draw the following statistical data for our group $O(V)$:
\begin{corollary}
The order of $O(V)$ is $2^{n^2}\prod_{j=1}^n(2^{2j}-1)$ and $O(V)$ contains $2^{2n}-1$ transvections. There is exactly one transvection corresponding to each hyperplane in $V$ which contains the radical vector $e_0$, (i.e. $e_0$ is the non-zero vector in the polarization of $\xi_0$.)
\end{corollary}
\begin{proof}
This follows directly from the Benson--Crawley-Boevey--Neeman theorem. Note that for a vector space over $\F_2$, the only possible reflections are transvections, and it is easy to see that any transvection in $O(V)$ must fix the radical vector $e_0$. Since only transvections are involved, the ramification formula for $r$ in the Hilbert series simplifies to
$$r=\sum_W\alpha_W$$ where $W$ runs through the hyperplanes of $V$ and $\alpha_W=\log_2|G_W|$ where $G_W$ is the pointwise stabiliser of $W$. Since $O(V)$ acts transitively on the set of hyperplanes containing $e_0$, a simple counting argument tells us that there is exactly one transvection associated to each.
\end{proof}

We turn next to the ring of invariants $S(U^*)^{O^-}$ for an orthogonal group $O^-$ which is the group of automorphisms of a $2n$-dimensional vector space $U$ endowed with a quadratic form 
$\xi_-$ of $-$type (Arf invariant 1).

\begin{theorem}
Let $T^\dagger$ be an abstract polynomial ring on generators $$\xi_0,\dots,\xi_{2n-2},d_{2n-1},\dots,d_{n+1}$$ where $\xi_i$ has degree $2^i+1$ and $d_j$ has degree $2^{2n-1}-2^{j-1}$. Define the degree preserving map
$$T^\dagger\to S(U^*)$$ by sending $\xi_0$ to the quadratic form $\xi_-$, sending $\xi_i$ to $Sq^{2^{i-1}}Sq^{2^{i-2}}\dots Sq^1\xi_0$ and sending $d_j$ to the symmetric polynomial of the same degree in the set of vectors in $V^*$ of $-$type (i.e. the elements $x$ of $V^*$ such that $\xi_-+x^2$ has $-$type). Then the image of $T^\dagger$ in $S(U^*)$ is equal to the ring of invariants $S(U^*)^{O^-}$  of the orthogonal group of automorphisms preserving $\xi_-$. The kernel of the map is generated by a regular sequence of $n-2$ elements which are homogeneous of degrees $2^{2n-1}+2^j$ for $1\le j\le n-2$. 

The $n-2$ relations can be expressed in matrix form as the sum of three column vectors. When they are expressed this way, we use an additional (redundant) generator $\xi_{2n-1}$ and impose an additional relation at the beginning which amounts to an expression for $\xi_{2n-1}$ in terms of the chosen generators:
$$\mathbf{T}_{n-1}
\left(\begin{matrix}d_{n+1}\\\vdots\\d_{2n-1}\end{matrix}\right)+
\left(\begin{matrix}\xi_{2n-1}\\\vdots\\\xi_{n+1}^{2^{n-2}}\end{matrix}\right)+
\sqrt{\left(\mathbf{L}'_n\mathbf{E}_n+\left(\mathbf{L}'_n\mathbf{K}_n+
\mathbf{R}'_n\right)\mathbf{F}_n\right)}.$$ 
The matrix 
$\mathbf{T}_{n-1}$ involves only the $\xi$'s and has determinant equal to $\Omega^-_{2n-2}(\xi_-)$, a certain polynomial which is intimately related to the sets of quadratic forms of 
$-$type on spaces of dimensions $2n-2$ and $2n$.
The second of the three column vectors has $i$th entry $\xi_{2n-i}^{2^{i-1}}$ and the $i$th relation can be interpreted as saying that in the ring $S$, $\xi_{2n-i}^{2^{i-1}}$ can be expressed in terms of other generators and lower powers of $\xi_{2n-i}$. This relation cannot be deduced from any of the other relations. 

The regular sequence of relations in $T^\dagger$ can be extended to a regular sequence of length equal to $3n-1$, the Krull dimension of $T^\dagger$, by taking the further $2n+1$ elements
$$\xi_0,\dots,\xi_n,d_{2n-1},\dots,d_{n+1}.$$
The ring of invariants $S(U^*)^{O^-}$ is a complete intersection.
\end{theorem}

We can read off corollaries about the Hilbert series as before.

\begin{corollary}
The Hilbert series of the ring $S(U^*)^{O^-}$ is
$$\frac{(1-t^{2^{2n-1}+2})\dots(1-t^{2^{2n-1}+2^{n-2}})}
{(1-t^2)(1-t^3)\dots(1-t^{2^{2n-2}+1})
\cdot(1-t^{2^{2n-1}-2^{2n-2}})\dots(1-t^{2^{2n-1}-2^{n}})}.$$
\end{corollary}
\begin{corollary}
The order of $O^-$ is $2^{n^2-n+1}(2^n+1)\prod_{j=1}^{n-1}(2^{2j}-1)$ and $O^-$ contains $2^{2n-1}+2^{n-1}$ transvections. There is exactly one transvection corresponding to each hyperplane in $U$ which is the kernel of a $+$type vector $x\in U^*$ (i.e. $\xi_-+x^2$ has $+$type.)
\end{corollary}

Finally we have the groups of $+$type. Let $U$ be a $2n$-dimensional vector space endowed with a quadratic form $\xi_+$ of $+$type (Arf invariant $0$). Notice that we still use vectors of $-$type to describe the generators $d_j$ even though this is the $+$type case!
\begin{theorem}
Let $T^\dagger$ be an abstract polynomial ring on generators $$\xi_0,\dots,\xi_{2n-2},d_{2n-1},\dots,d_{n}$$ where $\xi_i$ has degree $2^i+1$ and $d_j$ has degree $2^{2n-1}-2^{j-1}$. Define the degree preserving map
$$T^\dagger\to S(U^*)$$ by sending $\xi_0$ to the quadratic form $\xi_+$, sending $\xi_i$ to $Sq^{2^{i-1}}Sq^{2^{i-2}}\dots Sq^1\xi_+$ and sending $d_j$ to the symmetric polynomial of the same degree in the set of vectors in $V^*$ of $-$type (i.e. the elements $x$ of $V^*$ such that $\xi_++x^2$ has $-$type). Then the image of $T^\dagger$ in $S(U^*)$ is equal to the ring of invariants $S(U^*)^{O^+}$  of the orthogonal group of automorphisms preserving $\xi_+$. The kernel of the map is generated by a regular sequence of $n-1$ elements which are homogeneous of degrees $2^{2n-1}+2^j$ for $1\le j\le n-1$. 

The $n-1$ relations can be expressed in matrix form as the sum of three column vectors. When they are expressed this way, we use an additional (redundant) generator $\xi_{2n-1}$ and impose an additional relation at the beginning which amounts to an expression for $\xi_{2n-1}$ in terms of the chosen generators:
$$\mathbf{M}_n
\left(\begin{matrix}d_n\\\vdots\\d_{2n-1}\end{matrix}\right)+\left(\begin{matrix}\xi_{2n-1}\\\vdots\\\xi_{n+1}^{2^{n-2}}\\\\\xi_n^{2^{n-1}}\end{matrix}\right)+
\left(\begin{matrix}\sqrt{\mathbf{L}'_n\mathbf{E}_n+\left(\mathbf{L}'_n\mathbf{K}_n+
\mathbf{R}'_n\right)\mathbf{F}_n}\\\\\xi_n^{2^{n-1}}+f_n\end{matrix}\right).$$
The matrix 
$\mathbf{M}_{n}$ involves only polynomials in the $\xi$'s and its determinant is $\Omega^+_{2n-2}(\xi_+)$, a certain polynomial which is intimately related to the sets of quadratic forms of $+$type on spaces of dimensions $2n-2$ and $2n$. The polynomial $f_n$ is also a polynomial involving only the $\xi$'s.
The second of the three column vectors has $i$th entry $\xi_{2n-i}^{2^{i-1}}$ and the $i$th relation can be interpreted as saying that in the ring $S$, $\xi_{2n-i}^{2^{i-1}}$ can be expressed in terms of other generators and lower powers of $\xi_{2n-i}$. This relation cannot be deduced from any of the other relations. 

The regular sequence of relations in $T^\dagger$ can be extended to a regular sequence of length equal to $3n-1$, the Krull dimension of $T^\dagger$, by taking the further $2n$ elements
$$\xi_+,\dots,\xi_{n-1},d_{2n-1},\dots,d_{n}.$$
The ring of invariants $S(U^*)^{O^+}$ is a complete intersection.
\end{theorem}

\begin{corollary}
The Hilbert series of the ring $S(U^*)^{O^+}$ is
$$\frac{(1-t^{2^{2n-1}+2})\dots(1-t^{2^{2n-1}+2^{n-1}})}
{(1-t^2)(1-t^3)\dots(1-t^{2^{2n-2}+1})
\cdot(1-t^{2^{2n-1}-2^{2n-2}})\dots(1-t^{2^{2n-1}-2^{n-1}})}.$$
\end{corollary}
\begin{corollary}
The order of $O^+$ is $2^{n^2-n+1}(2^n-1)\prod_{j=1}^{n-1}(2^{2j}-1)$ and $O^+$ contains $2^{2n-1}-2^{n-1}$ transvections. There is exactly one transvection corresponding to each hyperplane in $U$ which is the kernel of a $-$type vector $x\in U^*$ (i.e. $\xi_++x^2$ has $-$type.)
\end{corollary}

\section{Connection with work of Domokos and Frenkel}

There is a potentially interesting connection between our calculations and those in \cite{df}.  If $\overline{\F_2}$ denotes the algebraic closure of $\F_2$ then we can consider the space $V\otimes\overline{\F_2}$ and its coordinate ring $\overline{\F_2}[V]\iso S(V^*)\otimes\overline{\F_2}$. Here we can look at the invariants of the full orthogonal group, the subgroup of $GL(V\otimes\overline{\F_2})$ preserving the quadratic form on $V$. Then the only invariant is the quadratic form itself. In this context, Domokos and Frenkel work with invariants of several vectors, that is, they study the invariants in the coordinate ring of a direct sum $V\otimes\overline{\F_2}\oplus\dots\oplus V\otimes\overline{\F_2}$ of several copies of $V$. This coordinate ring can be identified with the tensor product
$$\overline{\F_2}[V]\otimes\dots\otimes\overline{\F_2}[V]$$ 
and one can now consider the $\F_2$-linear map to $\overline{\F_2}[V]$ given by
$$s_1\otimes s_2\otimes s_3\otimes\dots\mapsto s_1s_2^2s_3^4\dots.$$
Using this it can be seen that the invariants of several vectors for the algebraic group give rise to the invariants $\xi_j$ for our finite group. This raises the possibility of extracting new information about invariants in our setting from the results of \cite{df}. For examle, Domokos and Frenkel show how to define an invariant which distinguishes the orthogonal group from the special orthogonal group and it seems reasonable to expect that this will map to an element of $S$ which distinguishes $O(V)$ from $SO(V)$: note that in characteristic $2$, the determinant  does not distinguish these groups. We do not carry through these investigations here. We thank Steve Donkin and Matias Domokos for drawing attention to this connection.

%%%%%%%%%%%%%%%%%%%%%%%%%%%%%%%%%%%%%%%%%%%%%%%%%%%%%%%%%%%%%%%%%%%%%%%%%%%%%%%%%%%%

\section{Invariants in the Symplectic Case}
\label{symplectic}

The invariant ring $S(U^*)^{Sp(U)}$ is known. This was first calculated by
Carlisle and Kropholler. Useful accounts have been published by Benson
(see Section 8.3 of \cite{benson}) and Neusel \cite{neusel}. We shall
have considerable need of this knowledge. 

\begin{proposition} \label{Carlisle Kropholler} The ring $S(U^*)^{Sp(U)}$
is generated by $\xi_1,\dots,\xi_{2n-1},c_{2n-1},\dots,c_n$. This ring
is a unique factorization domain. In terms of the stated generators, it
has a presentation given by a regular sequence $r_1,\dots,r_{n-1}$.
\end{proposition} 

To understand the relations, observe first that on multiplying both
sides of the matrix identity of Lemma \ref{prerelations} by $C_0B_0$ we
obtain the matrix identity 

$$\left( \begin{matrix} 0 & \xi_1 & \xi_2 & \xi_3 & \dots & \xi_{2n-1}
\\\\ \xi_1 & 0 & \xi_1^2 & \xi_2^2 & \dots & \xi_{2n-2}^2 \\\\ \xi_2 &
\xi_1^2 & 0& \xi_1^4 & \dots & \xi_{2n-3}^4 \\\\ \xi_3 & \xi_2^2 &
\xi_1^4 & 0 & \dots & \xi_{2n-4}^8 \\\\ \vdots &\vdots &\vdots &\vdots
&\ddots&\vdots\\\\ \xi_{2n-1} & \xi_{2n-2} ^2 & \xi_{2n-3} ^4 &
\xi_{2n-4} ^8 & \dots & 0\\ \end{matrix} \right)\left( \begin{matrix}
c_0\\\\ c_1\\\\ c_2\\\\ c_3\\\\ \vdots\\\\ c_{2n-1}\\ \end{matrix}
\right)=\left( \begin{matrix} \xi_{2n}\\\\ \xi_{2n-1}^2\\\\
\xi_{2n-2}^4\\\\ \xi_{2n-3}^8\\\\ \vdots\\\\ \xi_{1}^{2^{2n-1}}\\
\end{matrix} \right).$$ 

This matrix equation records $2n$ relations which hold in the ring
$S(U^*)^{Sp(U)}$. The first of these provides a formula for $\xi_{2n}$ in
terms of lower degree $\xi_i$ and the Dickson invariants, so telling us
that $\xi_{2n}$ may be omitted from the list of generators of
$S(U^*)^{Sp(U)}$. The next $n-1$ of these are, in disguise, the relations
$r_1,\dots,r_{n-1}$. More mysteriously, it turns out that the Dickson
invariants $c_{n-1},\dots,c_0$ of higher degree can all be expressed as
linear combinations of the Dickson invariants $c_{2n-1},\dots,c_n$ of
lower degree with coefficients in the ring
$\F_2[\xi_1,\dots,\xi_{2n-1}]$.
\begin{lemma}\label{redundancy of cs}
There is an $n\times n$ matrix $\mathbf{K}_n$ such that
$$\left(\begin{matrix}c_0\\\vdots\\c_{n-1}\end{matrix}\right)=
\mathbf{K}_n\left(\begin{matrix}c_n\\\vdots\\c_{2n-1}\end{matrix}\right)+\mathbf{E}_n.$$
The entries of $\mathbf{K}_n$ and $\mathbf{E}_n$ are all expressible as polynomials in the $\xi$'s. The top entry of the column vector $\mathbf{E}_n$ is  a certain polynomial $\Lambda_{2n}$ in $\xi_1,\dots,\xi_{2n-1}$ which we review in the next section and which is equal to the Dickson invariant $c_0$ in the ring $S$. The matrix $\mathbf{K}_n$ has zeroes on and above the anti-diagonal. The matrix and column vector are related by the recursive block matrix formula
$$\mathbf{K}_n=\left(\begin{matrix}0\ \cdots\ 0&0\\\\\mathbf{K}_{n-1}^{*2}&\mathbf{E}_{n-1}^{*2}\end{matrix}\right)$$
where the notation $\mathbf{K}_{n-1}^{*2}$ denotes the matrix obtained by squaring every element of $\mathbf{K}_{n-1}$.
\end{lemma}
This was proved by induction, and the
proof is recorded unaltered in the accounts of Benson and Neusel. 
It
remains of some interest to acquire a conceptual insight into this
aspect of the ring $S(U^*)^{Sp(U)}$. 
The key relations in the symplectic case are the first $n$ equations from the matrix equation exhibited following Proposition \ref{Carlisle Kropholler}. 
Partition the matrix into four $n\times n$ blocks. We are only concerned with the two blocks at the top which we denote by $\mathbf{L}_n$ and $\mathbf{R}_n$. The relations can now be expressed in matrix form as
$$\mathbf{L}_n\left(\begin{matrix}c_0\\\vdots\\c_{n-1}\end{matrix}\right)+
\mathbf{R}_n\left(\begin{matrix}c_n\\\vdots\\c_{2n-1}\end{matrix}\right)=
\left(\begin{matrix}\xi_{2n}\\\vdots\\\xi_{n+1}^{2^{n-1}}\end{matrix}\right)$$
We now make the substitutions for the redundant Dickson invariants $c_{n-1},\dots,c_0$ using Lemma \ref{redundancy of cs}. We then have
\begin{statement}\label{frs}{\bf The fundamental relations for the symplectic invariants $S(U^*)^{Sp(U)}$:}\end{statement}
$$\left(\mathbf{L}_n\mathbf{K}_n+\mathbf{R}_n\right)
\left(\begin{matrix}c_n\\\vdots\\c_{2n-1}\end{matrix}\right)=
\left(\begin{matrix}\xi_{2n}\\\vdots\\\xi_{n+1}^{2^{n-1}}\end{matrix}\right)+\mathbf{L}_n\mathbf{E}_n.$$
Here, the first relation simply gives expression for $\xi_{2n}$ in terms of other generators, so we can omit this relation and discard the redundant generator $\xi_{2n}$. The remaining $n-1$ relations are the relations $r_1,\dots,r_{n-1}$ referred to in 
Proposition \ref{Carlisle Kropholler}.

From here it is easy to establish the ring of invariants of $Sp(V)$. 

\begin{lemma} \label{eta} The ring of invariants of $Sp(V)$ is generated
by $S(U^*)^{Sp(U)}$ together with the single additional element
$$\eta:=\prod_{x\in V^*\setminus U^*}x.$$ Abstractly this is a
polynomial ring in one variable of degree $2^{2n}$ over $S(U^*)^{Sp(U)}$.
\end{lemma} \begin{proof} As we remarked following Definition
\ref{groups} the natural surjection $$Sp(V)\to Sp(U)$$ has kernel an
elementary abelian $2$-group $E$ of transvections. This subgroup acts
trivially on $U^*$ and clearly also fixes $\eta$. Hence
$$S^E\supseteq\F_2[x_1,\dots,x_{2n},\eta],$$ and a simple Galois
theoretic argument shows that equality holds. Now, the action of $Sp(V)$
on $S$ induces an action of $Sp(U)$ on $S^E$. The new element $\eta$ is
fixed by $Sp(U)$ and the action of $Sp(U)$ on the polynomials in
$x_1,\dots,x_{2n}$ is simply the classical action studied by Carlisle
and Kropholler. Hence the result follows. \end{proof} 

\begin{lemma}\label{6.3} The elements $\xi_0,\dots,\xi_{2n}$ are algebraically
independent. \end{lemma} \begin{proof} Since $\xi_1,\dots,\xi_{2n}$ all
belong to $S(U^*)$ while $\xi_0$ involves the additional variable $x_0$,
we need only show that $\xi_1,\dots,\xi_{2n}$ are algebraically
independent elements of $S(U^*)$. The determinant of the Jacobian matrix
$\left(\frac{\partial\xi_i}{\partial x_j}\right)_{1\le i,j\le2n}$ is
$\det C_0=c_0\not=0$ and the result follows from Proposition 5.4.2 of
\cite{benson}. \end{proof} 

By contrast, the elements $\xi_0,\dots,\xi_{2n+1}$ are obviously not
algebraically independent since they live in the ring $S$ of Krull
dimension $2n+1$ and they are $2n+2$ in number. At the risk of causing
untold confusion we shall bravely work on the assumption that the ring
$$\F_2[\xi_0,\xi_1,\xi_2,\dots]$$ really is an abstract polynomial ring
in the stated generators. The reason why we can get away with this
apparent travesty is that in any given situation we shall only be
concerned with the $\xi$ up to $\xi_{2n}$. On the other hand we shall be
proving our results in many cases by induction on $n$. 

%%%%%%%%%%%%%%%%%%%%%%%%%%%%%%%%%%%%%%%%%%%%%%%%%%%%%%%%%%%%%%%%%%%%%%%%
%%%%%%%%%%% 

\section{Some families of polynomials arising from determinants} 

Let $m$ be a positive integer. In the abstract commutative polynomial
ring $$\Z[X,\xi_1,\xi_2,\xi_3,\dots],$$ consider the polynomial 

$$H_m=\left| \begin{matrix} 2X & \xi_1 & \xi_2 & \xi_3 & \dots & \xi_{m}
\\\\ \xi_1 & 2X^2 & \xi_1^2 & \xi_2^2 & \dots & \xi_{m-1}^2 \\\\ \xi_2 &
\xi_1^2 & 2X^4 & \xi_1^4 & \dots & \xi_{m-2}^4 \\\\ \xi_3 & \xi_2^2 &
\xi_1^4 & 2X^8 & \dots & \xi_{m-3}^8 \\\\ \vdots &\vdots &\vdots &\vdots
&\ddots&\vdots\\\\ \xi_{m} & \xi_{m-1} ^2 & \xi_{m-2} ^4 & \xi_{m-3} ^8
& \dots & 2X^{2^{m}} \\ \end{matrix} \right|$$ 

This is the determinant of a symmetric matrix. On passing to the
quotient ring $$\F_2[X,\xi_1,\xi_2,\xi_3,\dots],$$ the matrix is
alternating. Since alternating matrices have even rank, it follows that
the determinant is zero modulo $2$ whenever $m$ is even, and we can make
the following definition: 

\begin{definition}\label{om def} For each even integer $m\ge0$, we write $\Omega_{m}(X)$ for
the image of the polynomial $\frac12H_m$ in
$\F_2[X,\xi_1,\xi_2,\xi_3,\dots].$ \end{definition} 
This observation has also been made by Domokos and Frenkel, see Proposition 4.11 of \cite{df}.
As an example, in case $m=2$ we find that
$$\Omega_2(X)=\xi_1^2X^4+\xi_2^2X^2+\xi_1^4X+\xi_1^3\xi_2.$$

When $m$ is odd, the image of $H_m$ in $\F_2[X,\xi_1,\xi_2,\xi_3,\dots]$
is non-zero and does not involve $X$. In fact it is the square of a
polynomial in $\F_2[\xi_1,\xi_2,\xi_3,\dots]$. The determinant of any alternating matrix over a commutative ring is a square; namely the square of the Pfaffian. The determinant of an alternating matrix over a commutative $\F_2$-algebra is more obviously square because the only contributing terms come from diagonally symmetric choices of elements from the matrix.
So we make the definition 
\begin{definition} For each even integer $m$, we write $\Lambda_{m}$ for
the square root of the image of the polynomial $H_{m-1}$ in
$\F_2[\xi_1,\xi_2,\xi_3,\dots],$ that is, the Pfaffian of the matrix defining $H_{m-1}$. \end{definition} 
For example, \begin{eqnarray*} \Lambda_2&=&\xi_1\\
\Lambda_4&=&\xi_1^5+\xi_2^3+\xi_1^2\xi_3\\
\Lambda_6&=&\xi_5\xi_3^{2}\xi_1^{4} +\xi_5\xi_2^{6} +\xi_5\xi_1^{10}
+\xi_4^{3}\xi_1^{4} +\xi_4^{2}\xi_3\xi_2^{4} \\
&&+\xi_4^{2}\xi_2\xi_1^{8} +\xi_4\xi_3^{4}\xi_2^{2}
+\xi_4\xi_2^{8}\xi_1^{2} +\xi_3^{7} +\xi_3^{4}\xi_1^{9} \\
&&+\xi_3^{2}\xi_2^{9} +\xi_3\xi_1^{18} +\xi_2^{12}\xi_1
+\xi_2^{3}\xi_1^{16} +\xi_1^{21} \\ \end{eqnarray*} 
In general, $\Lambda_{2n}$ has $\frac{(2n)!}{2^nn!}$ terms, this being the number of permutations in the symmetric group on $\{1,2,\dots,2n\}$ which are products of $n$ disjoint transpositions. If 
$$(i_1\ i_2)(i_3\ i_4)\dots(i_{2n-1}\ i_{2n})$$ is such a permutation then there is a corresponding contribution
$$\xi_{|i_1-i_2|}^{2^{\min\{i_1,i_2\}-1}}\xi_{|i_3-i_4|}^{2^{\min\{i_3,i_4\}-1}}\dots\xi_{|i_{2n-1}-i_{2n}|}^{2^{\min\{i_{2n-1},i_{2n}\}-1}}.$$
For example, $\Lambda_8$ has $105$ terms of which the leading term (giving $\xi_7$ the highest priority and $\xi_1$ the lowest) $\xi_7\xi_5^2\xi_3^4\xi_1^8$ arises from the contribution of the permutation 
$$(1\ 8)(2\ 7)(3\ 6)(4\ 5).$$

%%%%%%%%%%%%%%%%%%%%%%%%%%%%%%%%%%%%%%%%%%%%%%%%%%%%%%%%%%%%%%%%%%%%%%%%
%%%%%%%%%%% 

\section{How to understand $\Lambda_{m}$} 

Working in $S$, recall that the matrix $C_0$, defined in Section \ref{section definitions},
has determinant equal to
the Dickson invariant $c_0$. 

The extended matrix $C$ 
delivers a
sequence of square matrices $D_0,D_1\dots,D_{2n-1},D_{2n}=C_0$ where
$D_i$ is obtained by omitting the $i$th row of $C$. It is known that $D_i$ has determinant $c_0c_i$.

As in Section \ref{section definitions}, $B_0$ denotes the $2n\times2n$ matrix with
$(i,j)$-entry $b(e_i,e_j)$, $i,j\ge1$. Then $B_0$ is a non-singular
alternating matrix and so it has determinant $1$. Moreover
$$C_0^TB_0C_0=\left( \begin{matrix} 0 & \xi_1 & \xi_2 & \xi_3 & \dots &
\xi_{2n-1} \\\\ \xi_1 & 0 & \xi_1^2 & \xi_2^2 & \dots & \xi_{2n-2} ^2
\\\\ \xi_2 & \xi_1^2 &0 & \xi_1^4 & \dots & \xi_{2n-3} ^4 \\\\ \xi_3 &
\xi_2^2 & \xi_1^4 & 0 & \dots & \xi_{2n-4} ^8 \\\\
\vdots&\vdots&\vdots&\vdots& \ddots & \vdots \\\\ \xi_{2n-1} &
\xi_{2n-2}^2 & \xi_{2n-3}^4 & \xi_{2n-4}^8 & \dots& 0 \\\\ \end{matrix}
\right)$$ 

On taking determinants, noting that $\det B_0=1$, we find that
$$\det(C_0^TB_0C_0)=c_0^2$$ can be expressed as a polynomial in the
$Sp(U)$-invariants $\xi_1,\dots,\xi_{2n-1}$. As this matrix is clearly
congruent to the matrix of $H_{2n}$ modulo 2, it follows that $\Lambda
_{2n}=c_{0}$ in $S$ and that $c_{0}$ itself can be expressed in terms of
the $\xi _{i}$. 

Playing this game with $C$ in place of $C_0$, we have $$C^TB_0C=\left(
\begin{matrix} 0 & \xi_1 & \xi_2 & \xi_3 & \dots & \xi_{2n-1} & \xi_{2n}
\\\\ \xi_1 & 0 & \xi_1^2 & \xi_2^2 & \dots & \xi_{2n-2} ^2 & \xi_{2n-1}
^2 \\\\ \xi_2 & \xi_1^2 &0 & \xi_1^4 & \dots & \xi_{2n-3} ^4 &
\xi_{2n-2} ^4 \\\\ \xi_3 & \xi_2^2 & \xi_1^4 & 0 & \dots & \xi_{2n-4} ^8
& \xi_{2n-3} ^8 \\\\ \vdots&\vdots&\vdots&\vdots& \ddots & \vdots
&\vdots \\\\ \xi_{2n-1} & \xi_{2n-2}^2 & \xi_{2n-3}^4 & \xi_{2n-4}^8 &
\dots& 0& \xi_1^{2^{2n}}\\\\ \xi_{2n} & \xi_{2n-1}^2 & \xi_{2n-2}^4 &
\xi_{2n-3}^8 & \dots& \xi_1^{2^{2n}}&0 \\\\ \end{matrix} \right)$$ 

\begin{definition} We define polynomials $\Lambda_{2n,i}$ for each
$n\ge2$ and $0\le i\le2n$ by 
$$\Lambda _{2n,i} =Sq^{2^{2n}-2^{i}}\left( \Lambda _{2n}\right).$$
\end{definition} 

\begin{lemma} \label{lambdas}
\begin{enumerate} 
\item For each $i$ in the range $0\le i\le 2n$, we have $\Lambda _{2n,i}=c_{0}c_i$. Each $\Lambda_{2n,i}$ can also be interpreted as Pfaffians coming from the appropriate $2n\times 2n$ matrix obtained by omitting a row and corresponding column from $C^TB_0C$.
\item $\Lambda_{2n}$ belongs to the ring $\F_2[\xi_1,\dots,\xi_{2n-1}]$, and here it is irreducible and also linear in $\xi_{2n-1}$:
$$\Lambda_{2n}=\xi_{2n-1}\left(\Lambda_{2n-2}\right)^2+\text{ terms involving }\xi_1,\dots,\xi_{2n-2}.$$ Moreover, $$\Lambda_{2n,2n}=\Lambda_{2n},$$ and 
$$\Lambda_{2n,0}=\left(\Lambda_{2n}\right)^2.$$

\item For $1\le i\le 2n-1$, the polynomial $\Lambda_{2n,i}$ belongs to the ring $\F_2[\xi_1,\dots,\xi_{2n}]$ and is linear in $\xi_{2n}$:
$$\Lambda _{2n,i}=\xi_{2n}\left(\Lambda_{2n-2,i-1}\right)^2+\text{ terms involving }\xi_1,\dots,\xi_{2n-1}.$$ Moreover, $\Lambda_{2n,i}$ is not divisible by $\Lambda_{2n}$ for these values of $i$.

%\item The following identities hold in 
%$\mathbb{F}_{2}[\xi_1,\dots,\xi_{2n}]$: 
%\begin{eqnarray*} \Lambda _{2n}&=&\xi _{2n-1}\Lambda
%_{2n-2}^2+\xi _{2n-2}\Lambda_{2n-2,2n-3}^2+\dots+\xi_2\Lambda_{2n-2,1}^2+\xi_1\Lambda_{2n-2,0}^2 \\
%\Lambda_{2n,2n-1}&=&\xi_{2n}\Lambda_{2n-2,0}+\xi_{2n-1}^2\Lambda_{2n-2,1}+\xi_{2n-2}^4\Lambda_{2n-2,2}+\dots+\xi_2^{2^{2n-2}}\Lambda_{2n-2,2n-2}\\
%\Lambda_{2n,1}&=&\xi_{2n}\Lambda_{2n-2,2n-2}^4+\xi_{2n-1}\Lambda_{2n-2,2n-3}^4+\xi_{2n-2}\Lambda_{2n-2,2n-4}^4+\dots+\xi_2\Lambda_{2n-2,0}^4.
%\end{eqnarray*}
\end{enumerate} \end{lemma} 

\begin{proof}
\begin{enumerate} 
\item We noted above that 
$\Lambda_{2n}=c_{0}$. By Lemma \ref{steenroddickson} it follows that 
$\Lambda_{2n,i}=c_{0}c_{i}$. 
\item The nature of $\Lambda_{2n}$ as a polynomial in $\xi_1,\dots,\xi_{2n-1}$ is easily deduced from its definition as the square root of a determinant. Since the symplectic group acts transitively on the vectors of $U^*\setminus\{0\}$, it follows that in the invariant ring, $\Lambda_{2n}=c_0$ is irreducible and hence $\Lambda_{2n}$ is irreducible when viewed as a polynomial in the symplectic invariants $\xi_i$. 
\item From Lemma \ref{squarexi} and its Corollary we know that when a Steenrod operation is applied to $\Lambda_{2n}$ the only way in which $\xi_{2n}$ can become involved is through the application of $Sq^{2^{2n-1}}$ to $\xi_{2n-1}$, and taking this together with the Cartan formula, we calculate
\begin{eqnarray*}
\Lambda_{2n,i}&=&Sq^{2^{2n}-2^i}\Lambda_{2n}\\
&=&Sq^{2^{2n}-2^i}\left(\xi_{2n-1}\Lambda_{2n-2,i-1}^2+
\text{ terms involving }\xi_1,\dots,\xi_{2n-2}\right)\\
&=&Sq^{2^{2n-1}}\xi_{2n-1}\cdot Sq^{2^{2n-1}-2^i}\left(\Lambda_{2n-2,i-1}^2\right)+
\text{ terms involving }\xi_1,\dots,\xi_{2n-1}\\
&=&\xi_{2n}\cdot\left(Sq^{2^{2n-2}-2^{i-1}}\Lambda_{2n-2,i-1}\right)^2+
\text{ terms involving }\xi_1,\dots,\xi_{2n-1}\\
&=&\xi_{2n}\Lambda_{2n-2,i-1}^2+\text{ terms involving }\xi_1,\dots,\xi_{2n-1}.
\end{eqnarray*}
The last remarks now follow easily.

%\item Look at $\det C_{0}^{T}B_{0}C_{0}=\Lambda _{2n}$ and expand it
%along the top row. It is symmetric, so each term $\xi _{i}\xi _{j}f,$
%where $\xi _{i}$ is from the top row and $\xi _{j}$ from the first
%column, is cancelled by a term $\xi _{j}\xi _{i}f$ unless $i=j$. The
%coefficient of this term 
%$ \xi _{i}^{2}$ is the determinant of the
%matrix obtained from $C_{0}^{T}B_{0}C_{0}$ by deleting the 
%$\left( i+1\right)$st row and column, which is the matrix 
%$\left(M_{2n-2,i-1}\right) ^{2}$ with determinant 
%$ \left( \Lambda_{2n-2,i-1}\right) ^4.$ As squaring is additive over 
%$\mathbb{F}_{2}$, the first result follows. The other two identities can be proved in a similar way.
\end{enumerate}%
\end{proof} 
Examples of the $\Lambda_{2n,i}$ are given in Figure \ref{fig1}.

\begin{lemma} \label{J}
Let $J=\{s\in\F_2[\xi_0,\dots,\xi_{2n}];\
sc_i\in\F_2[\xi_0,\dots,\xi_{2n}]$ for each $i\}$. Then $J$ is the
principal ideal of $\F_2[\xi_0,\dots,\xi_{2n}]$ generated by
$\Lambda_{2n}$. \end{lemma} 

\begin{proof} Let $t\in J.$ Then for any $i$, we have $t\cdot\Lambda
_{2n.i}=t\cdot c_{0}c_{i}=tc_{i}\cdot\Lambda _{2n}$ and this is an
equation in\emph{\ }$\mathbb{F}_{2}[\xi _{0},\dots ,\xi _{2n}].$ Since $\Lambda _{2n}$ does not
divide $\Lambda _{2n,i}$, it must divide
$t.$ \end{proof} 

\begin{lemma} \label{squares} 
\begin{enumerate} \item If $s\in S$ has
the property that $s^2$ can be expressed as a polynomial in at most $2n$
of $\xi_0,\dots,\xi_{2n}$, then $s$ itself is a polynomial generated by
the same $\xi_i$. \item If $f$ is a polynomial of degree $2^{2n+1}$ in
$F_2[\xi_0,\dots,\xi_{2n}]$ which is a square in the ambient ring $S$
and which involves $\xi_{2n}$ then $\frac{\partial f}{\partial
\xi_i}=c_0c_i$ for each $i$ in the range $0$ to $2n$. \end{enumerate}
\end{lemma} 

\begin{proof} 
\begin{enumerate} \item An element $f$ of $S$ is a square
if and only if $$\frac{\partial f}{\partial x_i}=0$$ for each $i$ in the
range $0$ to $2n$. In the light of the first part of this lemma, we also
know that if $f$ belongs to the subring $\F_2[\xi_0,\dots,\xi_{2n}]$
then $f$ is intrinsically a square within this ring if and only if
$$\frac{\partial f}{\partial \xi_i}=0$$ for each $i$ in the range $0$ to
$2n$. The transition between these two conditions is made via the
Jacobian identity

$$\left( \begin{matrix} \tfrac{\partial f}{\partial x_0} \\\\
\tfrac{\partial f}{\partial x_1} \\\\ \tfrac{\partial f}{\partial x_2}
\\\\ \tfrac{\partial f}{\partial x_3} \\\\ \vdots\\\\ \tfrac{\partial
f}{\partial x_{2n}} \\ \end{matrix} \right)=\left( \begin{matrix}
\tfrac{\partial\xi_0}{\partial x_0} & \tfrac{\partial\xi_1}{\partial
x_0} & \tfrac{\partial\xi_2}{\partial x_0} &
\tfrac{\partial\xi_3}{\partial x_0} & \dots &
\tfrac{\partial\xi_{2n}}{\partial x_0} \\\\
\tfrac{\partial\xi_0}{\partial x_1} & \tfrac{\partial\xi_1}{\partial
x_1} & \tfrac{\partial\xi_2}{\partial x_1} &
\tfrac{\partial\xi_3}{\partial x_1} & \dots &
\tfrac{\partial\xi_{2n}}{\partial x_1} \\\\
\tfrac{\partial\xi_0}{\partial x_2} & \tfrac{\partial\xi_1}{\partial
x_2} & \tfrac{\partial\xi_2}{\partial x_2} &
\tfrac{\partial\xi_3}{\partial x_2} & \dots &
\tfrac{\partial\xi_{2n}}{\partial x_2} \\\\
\tfrac{\partial\xi_0}{\partial x_3} & \tfrac{\partial\xi_1}{\partial
x_3} & \tfrac{\partial\xi_2}{\partial x_3} &
\tfrac{\partial\xi_3}{\partial x_3} & \dots &
\tfrac{\partial\xi_{2n}}{\partial x_3} \\\\ \vdots &\vdots &\vdots
&\vdots &\ddots&\vdots\\\\ \tfrac{\partial\xi_0}{\partial x_{2n}} &
\tfrac{\partial\xi_1}{\partial x_{2n}} & \tfrac{\partial\xi_2}{\partial
x_{2n}} & \tfrac{\partial\xi_3}{\partial x_{2n}} & \dots &
\tfrac{\partial\xi_{2n}}{\partial x_{2n}} \\ \end{matrix} \right)\left(
\begin{matrix} \tfrac{\partial f}{\partial \xi_0} \\\\ \tfrac{\partial
f}{\partial \xi_1} \\\\ \tfrac{\partial f}{\partial \xi_2} \\\\
\tfrac{\partial f}{\partial \xi_3} \\\\ \vdots\\\\ \tfrac{\partial
f}{\partial \xi_{2n}} \\ \end{matrix} \right).$$ 

The Jacobian matrix is equal to $$B\widehat C^T.$$ Hence, if $f$ is a
square in $S$ then the Jacobian identity tells us that for each $i$ in
the range $1$ to $2n$, $$ \tfrac{\partial f}{\partial
\xi_{2n}}x_i^{2^{2n}}+ \tfrac{\partial f}{\partial
\xi_{2n-1}}x_i^{2^{2n-1}}+\dots+ \tfrac{\partial f}{\partial
\xi_{1}}x_i^{2}+ \tfrac{\partial f}{\partial \xi_{0}}x_i=0.$$ From Lemma
\ref{dick} we deduce that $$ \tfrac{\partial f}{\partial \xi_{i}}=
\tfrac{\partial f}{\partial \xi_{2n}}c_i$$ for each $i\ge0$. The
hypotheses here guarantee that there is at least one choice of $i$ for
which $\tfrac{\partial f}{\partial \xi_{i}}=0$ and the above equations
now show that all $\tfrac{\partial f}{\partial \xi_{i}}$ vanish. Thus
$f$ is an intrinsic square as required. 
\item  
From the proof of the last part, we know that $\tfrac{\partial f}{%
\partial \xi _{i}}=\tfrac{\partial f}{\partial \xi _{2n}}c_{i}$. Multiplying
both sides of this equation by $c_{0}=\Lambda _{2n}$ in $S,$ recalling that $%
\Lambda _{2n,i}=c_{0}c_{i},$  gives an equality in $F_{2}[\xi _{0},\dots
,\xi _{2n}]$: 
\begin{equation*}
\Lambda _{2n}\tfrac{\partial f}{\partial \xi _{i}}=\tfrac{\partial f}{%
\partial \xi _{2n}}\Lambda _{2n,i}.
\end{equation*}
Since $f$ involves $\xi_{2n}$ and has degree less than that of $\xi_{2n}^2$ we know that it is linear in $\xi_{2n}$ and that 
$\tfrac{\partial f}{\partial \xi _{2n}}$ is the coefficient of $\xi_{2n}$. Since $\Lambda_{2n}$ does not divide $\Lambda_{2n,i}$ unless $i=0$ or $2n$ we deduce that $\Lambda_{2n}$ divides $\tfrac{\partial f}{\partial \xi _{2n}}$ and on grounds of degree, the result follows.
\end{enumerate} \end{proof} 

%%%%%%%%%%%%%%%%%%%%%%%%%%%%%%%%%%%%%%%%%%%%%%%%%%%%%%%%%%%%%%%%%%%%%%%%
%%%%%%%%% 

%\newpage\ 

\begin{figure}%
\hrule
\begin{eqnarray*} \Lambda_2&=&\xi_1\\
\Lambda_{2,1}&=&\xi_2\\ \Lambda_{2,0}&=&\xi_1^2\\
\Lambda_4&=&\xi_3\xi_1^{2} +\xi_2^{3} +\xi_1^{5} \\
\Lambda_{4,3}&=&\xi_4\xi_1^{2} +\xi_3^{2}\xi_2 +\xi_2^{4}\xi_1 \\
\Lambda_{4,2}&=&\xi_4\xi_2^{2} +\xi_3^{3} +\xi_1^{9} \\
\Lambda_{4,1}&=&\xi_4\xi_1^{4} +\xi_3\xi_2^{4} +\xi_2\xi_1^{8}\\
\Lambda_{4,0}&=&\xi_3^{2}\xi_1^{4} + \xi_2^{6} + \xi_1^{10} \\
\Lambda_6&=&\xi_5\xi_3^{2}\xi_1^{4} +\xi_5\xi_2^{6} +\xi_5\xi_1^{10}
+\xi_4^{3}\xi_1^{4} +\xi_4^{2}\xi_3\xi_2^{4}+ \\
&&\quad\xi_4^{2}\xi_2\xi_1^{8} +\xi_4\xi_3^{4}\xi_2^{2}
+\xi_4\xi_2^{8}\xi_1^{2} +\xi_3^{7} +\xi_3^{4}\xi_1^{9}+ \\
&&\qquad\xi_3^{2}\xi_2^{9} +\xi_3\xi_1^{18} +\xi_2^{12}\xi_1
+\xi_2^{3}\xi_1^{16} +\xi_1^{21}\\
\Lambda_{6,5}&=&\xi_6\xi_3^{2}\xi_1^{4} +\xi_6\xi_2^{6} +\xi_6\xi_1^{10}
+\xi_5^{2}\xi_4\xi_1^{4}+\xi_5^{2}\xi_3\xi_2^{4}+\\ 
&&\quad\xi_5^{2}\xi_2\xi_1^{8} +\xi_4^{5}\xi_2^{2} +\xi_4^{4}\xi_3^{3}
+\xi_4^{4}\xi_1^{9} +\xi_4\xi_3^{8}\xi_1^{2}+\\ 
&&\qquad\xi_3^{10}\xi_2
+\xi_3^{8}\xi_2^{4}\xi_1+\xi_3\xi_2^{16}\xi_1^{2} +\xi_2^{19}
+\xi_2^{16}\xi_1^{5}\\ 
\Lambda_{6,4}&=&\xi_6\xi_4^{2}\xi_1^{4}
+\xi_6\xi_3^{4}\xi_2^{2} +\xi_6\xi_2^{8}\xi_1^{2}+\xi_5^{3}\xi_1^{4} +\xi_5^{2}\xi_3^{5}+\\
&&\quad\xi_5^{2}\xi_2^{9}
+\xi_5\xi_4^{4}\xi_2^{2}+\xi_5\xi_3^{8}\xi_1^{2} +\xi_4^{6}\xi_3
+\xi_4^{4}\xi_2^{8}\xi_1+\\
&&\qquad\xi_4^{2}\xi_3^{8}\xi_2+\xi_3^{12}\xi_1
+\xi_3\xi_1^{34} +\xi_2^{3}\xi_1^{32} +\xi_1^{37}\\
\Lambda_{6,3}&=&\xi_6\xi_4^{2}\xi_2^{4} +\xi_6\xi_3^{6} +\xi_6\xi_1^{18}
+\xi_5^{3}\xi_2^{4} +\xi_5^{2}\xi_4\xi_3^{4}+ \\
&&\quad\xi_5^{2}\xi_2\xi_1^{16} +\xi_5\xi_4^{4}\xi_3^{2}
+\xi_5\xi_2^{16}\xi_1^{2} +\xi_4^{7} +\xi_4^{4}\xi_1^{17}+ \\
&&\qquad\xi_4^{2}\xi_2^{17} +\xi_4\xi_1^{34} +\xi_3^{4}\xi_2^{16}\xi_1
+\xi_3^{2}\xi_2\xi_1^{32} +\xi_2^{4}\xi_1^{33}\\
\Lambda_{6,2}&=&\xi_6\xi_4^{2}\xi_1^{8} +\xi_6\xi_3^{2}\xi_2^{8}
+\xi_6\xi_2^{2}\xi_1^{16} +\xi_5^{3}\xi_1^{8} +\xi_5^{2}\xi_4\xi_2^{8}+\\
&&\quad\xi_5^{2}\xi_3\xi_1^{16} +\xi_5\xi_3^{10} +\xi_5\xi_2^{18}
+\xi_4^{3}\xi_3^{8} +\xi_4^{2}\xi_3\xi_2^{16}+\\
&&\qquad\xi_4\xi_2^{2}\xi_1^{32} +\xi_3^{8}\xi_1^{17} +\xi_3^{3}\xi_1^{32}
+\xi_2^{24}\xi_1 +\xi_1^{41}\\ 
\Lambda_{6,1}&=&\xi_6\xi_3^{4}\xi_1^{8}
+\xi_6\xi_2^{12} +\xi_6\xi_1^{20} +\xi_5\xi_4^{4}\xi_1^{8}
+\xi_5\xi_3^{8}\xi_2^{4}+ \\ 
&&\quad\xi_5\xi_2^{16}\xi_1^{4}
+\xi_4^{5}\xi_2^{8} +\xi_4^{4}\xi_3\xi_1^{16} +\xi_4\xi_3^{12}
+\xi_4\xi_1^{36}+\\ 
&&\qquad\xi_3^{8}\xi_2\xi_1^{16} +\xi_3^{5}\xi_2^{16}
+\xi_3\xi_2^{4}\xi_1^{32} +\xi_2^{25} +\xi_2\xi_1^{40} \\
\end{eqnarray*}
\hrule
\caption{Examples of the $\Lambda_{2n,i}$}
\label{fig1}
\end{figure}

\vfill

\ 
\newpage

%%%%%%%%%%%%%%%%%%%%%%%
%%%%%%%%%%%%%%%%%%%%%%%

\section{The Chern Polynomials} 

In this section we define Chern polynomials whose coefficients can be plainly seen to be invariants of the orthogonal group $O(V)$ and ``quadratic Chern polynomials'' whose coefficients are plainly invariants of the symplectic group $Sp(U)$. 

\begin{definition}\label{the ps}\  
Let $A^+$ denote the set  $$\{x\in V^*;\ \xi_0+x^2 \text{ has $+$type}\}$$ and let $A^-$ denote the set
$$\{x\in V^*;\ \xi_0+x^2 \text{ has $-$type}\}.$$  Let $A=A^+\cup A^-$. Note that $A=V^*\setminus U^*$. 
Polynomials
$P^+(t)$ and $P^-(t)$ in the polynomial ring $S[t]$ in one
variable $t$ of degree $1$ are defined as follows: 
\begin{eqnarray*}
\displaystyle P^+(t):=\prod_{x\in A^+}(t+x), \quad && \quad
P^-(t):=\prod_{x\in A^-}(t+x),\\ 
P(t)&:=&\prod_{x\in A}(t+x).
\end{eqnarray*}
We define $n$ particular invariants $d_{2n-1},\dots,d_n$ by picking certain coefficients of $P^-(t)$:
\begin{center}
{\em $d_j$ is the coefficient of $P^-(t)$ of degree $2^{2n-1}-2^{j-1}$}
\end{center}
for $j$ in the range $n\le j\le2n-1.$
\end{definition}

The orthogonal group $O(V)$ permutes the elements of $A^+$ and $A^-$, so it is clear that the coefficients of $P^+(t)$ and $P^-(t)$ belong to the invariant ring $S^{O(V)}$. Note also that $P(t)=P^+(t)P^-(t)$.

The quadratic Chern polynomials are closely related:
 \begin{definition}\label{the qs}
Let $B^+$ be the set of all quadratic forms on $V$ 
of $+$type and $B^-$ the set of all of $-$type. Let $B=B^+\cup B^-$. Note that the elements of $B$ belong to $S(U^*)$.  Note that, from the discussion in \S\ref{quadforms} we know that $B=\{\xi_0+x_0^2+x^2;\ x\in U^*\}$.
We define quadratic Chern polynomials  $Q^+(X)$, $Q^-(X)$ and $Q(X)$ in the polynomial ring $S[X]$ in one variable $X$ of degree $2$ as follows:  
\begin{eqnarray*}
Q^+(X):=\prod_{q\in B^+}(X+q), \quad && \quad
Q^-(X):=\prod_{q\in B^-}(X+q), \\
Q(X)&:=&\prod_{q\in B}(X+q).
\end{eqnarray*}
\end{definition}
The symplectic group $Sp(U)$ permutes the elements of $B^+$ and $B^-$. Thus the coefficients of $Q^+(X)$ and $Q^-(X)$  are invariants
of the symplectic group $Sp(V)$. Note also that $Q(X)=Q^+(X)Q^-(X)$.

We'll begin by illustrating these polynomials in the low dimensional cases:

\begin{example}The Case $n=1$ and $\dim V=3$. \end{example}

If $n=1$ then \begin{eqnarray*} Q^+(X)&=&X^3+c_1X^2+\xi_1^2,\\
Q^-(X)&=&X+c_1. \end{eqnarray*} In this case $c_1=x_1^2+x_1x_2+x_2^2$ happens to be the
unique quadratic from of $-$type, and we also have the relation
$c_0=\xi_1$. Thus $Q^+$ is also given by $$Q^+(X)=X^3+c_1X^2+c_0^2,$$
reflecting the coincidence $$Sp(U)=GL(U).$$ 
The invariant ring $S(U^*)^{Sp(U)}$ is generated by $\xi_1=c_0$ and $c_1$: it is the ring of Dickson invariants.

Further, with respect to a suitable basis, $\xi_0=x_0^2+x_1x_2$ and
\begin{eqnarray*}
P^+(t)&=&t^3 + (x_0+x_1+x_2)t^2 + \xi_0t + \xi_0(x_0+x_1+x_2)+\xi_1, \\
&=&t^3 + d_1t^2 + \xi_0t + \xi_0d_1+\xi_1, \\
P^-(t)&=&t+x_0+x_1+x_2,\\
&=&t+d_1.
\end{eqnarray*}

The ring $S^{O(V)}$ is a polynomial ring with generators $x_0+x_1+x_2$, $\xi_0$ and $\xi_1$. In fact $O(V)$ is isomorphic to the symmetric group on $3$ letters and its action on $V$ permutes the basis $e_1,e_2,e_0+e_1+e_2$. On $V^*$, it permutes the dual basis $x_0+x_1,x_0+x_2,x_0$, and the invariants $x_0+x_1+x_2$, $\xi_0$ and $\xi_1$ are the corresponding elementary symmetric polynomials.

\begin{example}\label{egn=2}The Case $n=2$ and $\dim V=5$. \end{example}

If $n=2$ then \begin{eqnarray*}
Q^+(X)&:=&X^{10}+\xi_1^2X^7+(c_3+\xi_1\xi_2)X^6+\xi_2^2X^5\\
&&\qquad\qquad+(c_2+\xi_1\xi_3+\xi_1^4)X^4+\xi_1^2\xi_2^2X^2+\xi_1^6X+
(\xi_1^4c_3+\xi_1^5\xi_2+\xi_2^4)\\
Q^-(X)&:=&X^6+\xi_1^2X^3+(c_3+\xi_1\xi_2)X^2+\xi_2^2X+(c_2+\xi_1\xi_3)\\
\end{eqnarray*} 
We have also computed the polynomials $P^-(t)$ and $P^+(t)$ in this case, at least in terms of the two coefficients $d_3,d_2$ of $P^-(t)$. We have
\begin{eqnarray*}
P^-(t)&=&t^6+\xi_0t^4+\xi_1t^3+d_3t^2+(\xi_2+\xi_1\xi_0)t+d_2,\\
P^+(t)&=&t^{10}+\xi_0t^8+\xi_1t^7+(d_3+\xi_0^2)t^6+(\xi_2+\xi_1\xi_0)t^5+(d_2+\xi_1^2+\xi_0^3)t^4+\\
&&\qquad \xi_1\xi_0^2t^3+(\xi_0^2d_3+\xi_2\xi_1)t^2+(\xi_2\xi_0^2+\xi_1\xi_0^3+\xi_1^3)t+\\
&&\qquad\qquad(\xi_0^2d_2+\xi_1^2d_3+\xi_2^2+\xi_2\xi_1\xi_0).
\end{eqnarray*}
In this case there are $6$ quadratic forms of minus type
and so $Q^-(X)$ has degree $12$. Since the ring of invariants of $Sp(U)$
is generated by $\xi_1,\xi_2,\xi_3,c_3,c_2$ subject to the single
relation
$$\xi_1^2c_2+\xi_2^2c_3+\xi_1^3\xi_3+\xi_1\xi_2^3+\xi_3^2+\xi_1^6=0,$$
we see that there is no ambiguity in expressing the coefficients of
$Q^-$ as polynomials in $\xi_1,\xi_2,\xi_3,c_3,c_2$. There are the following expressions for the other two Dickson invariants and for $\xi_4$ in terms of the minimal generating set.
\begin{eqnarray*}
c_1&=&\xi_1^2c_3+\xi_3\xi_2+\xi_2\xi_1^3,\\
c_0&=&\Lambda_4,\\
\xi_4&=&(\xi_3+\xi_1^3)c_3+\xi_2c_2+\xi_3\xi_2\xi_1+\xi_2\xi_1^4.
\end{eqnarray*}
When we introduce the orthogonal invariants $d_3$ and $d_2$ we find that
\begin{eqnarray*}
c_3&=&d_3^2+\xi_2\xi_1+\xi_1^2\xi_0+\xi_0^4\\
c_2&=&d_2^2+\xi_0^2d_3^2+\xi_3\xi_1+\xi_2^2\xi_0\\
\xi_3&=&\xi_1d_2+(\xi_2+\xi_1\xi_0)d_3+\xi_2\xi_0^2+\xi_1^3
\end{eqnarray*}
The orthogonal group $O(V)$ has invariant ring generated by $\xi_0,\xi_1,\xi_2,d_3,d_2$. Note that in this case, $O(V)$ is isomorphic to the symmetric group on $6$ letters and so admits an action on a $6$ dimensional space permuting a basis. We can then take $V$ to be the $5$ dimensional space consisting of the zero-sum vectors in the chosen basis and this gives the present representation of $O(V)$.

\begin{lemma}\label{p squared}  The following identities hold in $S[t]$:
\begin{eqnarray*} 
\left(P^+(t)\right)^2&=&Q^+(t^2+\xi_0)\\
\left(P^-(t)\right)^2&=&Q^-(t^2+\xi_0)\\
P(t)=P^-(t)P^+(t)&=&D(t+x_0)=D(t)+D(x_0)\\ 
\end{eqnarray*} 
\end{lemma} 

\begin{proof}
For the first two equalities, note that 
\begin{eqnarray*}
\left(P^{\pm}(t)\right)^2 &=& \left(\prod_{x\in A^{\pm}}(t+x)\right)^2 \\
&=& \prod_{x\in A^{\pm}}(t^2+x^2) \\
&=& \prod_{q\in B^{\pm}}(t^2+\xi_0 +q) \\
&=& Q^{\pm}(t^2+x_0) 
\end{eqnarray*} 
as $q\in B^{\pm}$ precisely when $q=\xi_0 +x^2$ with $x\in A^\pm$.

For the third, 
\begin{eqnarray*}
P^-(t)P^+(t) &=& \left(\prod_{x\in A^-}(t+x)\right)\left( \prod_{x\in A^+}(t+x)\right) \\
&=&\prod_{x\in A}(t+x) \\
&=& \prod_{x\in U^*}(t+x_0 +x) \\
&=& D(t+x_0). 
\end{eqnarray*} 
Further $D(X)$ is a polynomial in powers of 2, so is additive.
\end{proof}

\begin{lemma}\label{deg arg}
\begin{enumerate}
\item The coefficients of $Q^-(X)$ belong to the subring of $S$ generated by $$\xi_1,\dots,\xi_{2n-1},c_{2n-1},\dots c_n.$$ Moreover, they are linear in the Dickson invariants $c_{2n-1},\dots c_n.$
\item The polynomial $c_0Q^-(X)$ has all coefficients in the ring $\F_2[\xi_1,\dots,\xi_{2n}].$
\item The coefficients of $Q^+(X)$ belong to the subring of $S$ generated by $$\xi_1,\dots,\xi_{2n-1},c_{2n-1},\dots c_n.$$ Moreover, the only conceivable terms which are not linear in the Dickson invariants are terms involving $c_{2n-1}^2$.
\item The polynomial $c_0^2Q^+(X)$ has all coefficients in the ring $\F_2[\xi_1,\dots,\xi_{2n}].$
\item The squares of the coefficients of $P^-(t)$ and of $P^+(t)$ belong to the subring of $S$ generated by $\F_2[\xi_0,\xi_1,\dots,\xi_{2n-1},c_{2n-1},\dots,c_n].$

\end{enumerate}
\end{lemma}

\begin{proof}
\begin{enumerate}
\item The coefficients of $Q^-$ are symplectic invariants and so, using our knowledge of the invariant ring for that case, Proposition \ref{Carlisle Kropholler}, these coefficients lie in the subring $\F_2[\xi_1,\dots,\xi_{2n-1},c_{2n-1},\dots c_n]$. Since there are $2^{2n-1}-2^{n-1}$ quadratic forms of $-$type, the degree of $Q^-$ is $2(2^{2n-1}-2^{n-1})=2^{2n}-2^n$. On the other hand, the least degree of an element of $\F_2[\xi_1,\dots,\xi_{2n-1},c_{2n-1},\dots c_n]$ which is quadratic in the Dickson invariants is $\deg c_{2n-1}^2=2^{2n}$ and this is greater than the degree of $Q^-$. Hence the coefficients of $Q^-$ are at worst linear in the $c_j$.
\item We know that for each $j$, $c_0c_j$ belongs to $\F_2[\xi_1,\dots,\xi_{2n}]$ by Lemma \ref{lambdas}(i). Part (i) here says that the coefficients are linear in the $c_j$ and so the result follows.
\item As in (i) we can use Proposition \ref{Carlisle Kropholler}. For $n=1$ or $2$ we can see from the examples following Definition \ref{the qs} that the result holds. For $n\ge3$, the degree argument of (i) yields only the weaker stated result because there are more ($2^{2n-1}+2^{n-1}$) quadratic forms of $+$type. 
\item Again, Lemma \ref{lambdas}(i) applies, but for the moment we need the factor $c_0^2$ because of the weaker conclusion of (iii).
\item This follows from parts (i) and (iii) and Lemma \ref{p squared}.  
\end{enumerate}
\end{proof}

Note that $D(X)$ vanishes on $U^*$ and is constant on $V^*\setminus
U^*=A$. So the value $D(x_0)$ is not in fact dependent on $x_0$ and
equally, not dependent on any particular choice of $e_1,\dots,e_{2n}$.
In fact 
$$D(x_0)=\prod_A x$$ 
is the special
additional invariant $\eta$ of $Sp(V)$ introduced in Lemma \ref{eta}. Lemma \ref{p squared} shows that $$D(x_0)=P^+(0)P^-(0),$$ and so $D(x_0)$ also has
its square in the subring
$$\F_2[\xi_0,\xi_1,\dots,\xi_{2n-1},c_{2n-1},\dots,c_n].$$ 

We conclude this section with a remark about $Q^-$ which is needed later:
Let $W$ be a maximal $b$-isotropic subspace of $U$. Such a subspace of $U$ has dimension $n$. Any quadratic form polarizing to $b$ restricts to a linear functional on $W$ because its polarization vanishes on $W$. The quadratic forms of $-$type restrict to non-zero linear functionals on $W$ and every such linear functional on $W$ arises in this way. 
Since every automorphism of $W$ arises as the restriction of some symplectic automorphism of $V$ it follows that the restriction of $Q^-(t^2)$ to $W$ is a power of the Dickson polynomial for $W$. On grounds of degree we therefore have
\begin{lemma}\label{mis1}
The image of $Q^-(t^2)$ in $S(W^*)$ is 
$$\left(\prod_{0\ne x\in W^*}(t+x)\right)^{2^n}.$$
\end{lemma}

%%%%%%%%%%%%%%%%%%%%%%%%%%%%%%%%%%%%%%%%%%%%%%%%%%%%%%%%%%%%%%%%%%%%%%%

\section{How to understand $\Omega_m(X)$} 

We shall study the image of $\Omega_{m}(X)$ in the polynomial ring $S[X]$
over our symmetric algebra $S$, using the specialization
$$\F_2[X,\xi_1,\xi_2,\xi_3,\dots]\to S[X]$$ defined by $X\mapsto X$ and
$\xi_i\mapsto\xi_i$. 

\begin{lemma} \label{omegas} 
\begin{enumerate} 
\item 
$\Omega_{2n}(X)=\sum_{i=0}^{2n}\left(\Lambda_{2n,i}\right)^2 X^{2^i}+\delta$ 
 where $\delta \in  \F_2[\xi_1,\xi_2,\dots,\xi_{2n}]$. 
\item 
In the ring $S$ we have $\Omega_{2n}(X)=c_0^2Q(X).$
\item
$\Omega_{2n}(X)=c_0^2Q^-(X)Q^+(X)$, and $c_0$, $Q^-(X)$ and $Q^+(X)$ 
are irreducible elements of the
ring $S(U^*)^{Sp[U]}[X]$.
\item $c_0Q^-(X)$ and $c_0Q^+(X)$ both belong to $\F_2[X,\xi_1,\dots,\xi_{2n}]$.
\end{enumerate}
\end{lemma} 

\begin{proof} 
\begin{enumerate} 
\item Looking at the standard expansion
of the determinant $A_{2n}$ we see first that any term involving a product
of two or more of the diagonal entries will have a coefficient divisible
by $4$. So these make zero contribution to $\Omega_{2n}$. For  $0\le
i\le2^{2n}$, we see that the coefficient of $X^{2^i}$ in $\Omega_{2n}$ is
precisely the determinant of the matrix $A_{2n,i}$ obtained by omitting 
the $i$th row and column (counting from $0$ to $2n$) from $C^T B_0 C$. 
This determinant is equal to $(\Lambda_{2n,i})^2$ as noted in the proof
of Lemma \ref{lambdas}(iv).
\item 
Recall Definition \ref{the qs} that
$Q(X)=\prod_{q\in B}(X+q)$ where $B=\{\xi_0+x_0^2+x^2;\ x\in U^*\}.$
Using Dickson's Theorem (see Definition \ref{dickthm}), we know that the polynomial 
$D'(X):=\sum_{i=0}^{2n} c_i^2 X^{2^i}$ has zero set precisely 
$\{x^2;\ x\in U^*\}$. (Note that $D'(x^2)=D(x)^2$.) 
Thus 
$$Q(X)=D'(X+\xi_0+x_0^2)=D'(X)+D'(\xi_0+x_0^2).$$
We claim that $\Omega_{2n}=c_0^2 Q(X)$. First, it follows from Lemma \ref{lambdas}(i) that $c_0^2D'(X)=\sum_{i=0}^{2n}(\Lambda_{2n,i})^2X^{2^i}$ and by part (i), this coincides with the part of $\Omega_{2n}(X)$ which involves $X$. Therefore 
$$c_0^2Q(X)+\Omega_{2n}(X)$$ does not involve $X$ and to prove that it is zero it suffices to prove that
$$\Omega_{2n}(\xi_0+x_0^2)=0.$$
To this end we need to work over $\Z$ rather than $\F_2$ and we shall temporarily work with two abstract polynomial rings and the ring homomorphism as follows:
$$\alpha:\Z[X,\xi_1,\xi_2,\xi_3,\dots]\to\Z[x_1,x_2,\dots,x_{2n}]$$
where
$$\alpha(\xi_i)=\sum_{\ell=1}^n(x_{2\ell-1}^{2^i}x_{2\ell}+x_{2\ell-1}x_{2\ell}^{2^i}),$$
and
$$\alpha(X)=\sum_{\ell=1}^{n}x_{2\ell-1}x_{2\ell}.$$

Consider the matrices $C$ and $B_0C^T$ over $\Z$, and insert
respectively a row and
a column of zeros to make the matrices square. Then clearly they have
determinant equal to zero, and further we have the matrix equation:
$$
\left(
\begin{matrix} x_1&x_2&x_3&\dots&x_{2n}&0\\\\
x_1^2&x_2^2&x_3^2&\dots&x_{2n}^2&0\\\\
x_1^4&x_2^4&x_3^4&\dots&x_{2n}^4&0\\\\
\vdots&\vdots&\vdots&\ddots&\vdots&\vdots\\\\
x_1^{2^{2n}}&x_2^{2^{2n}}&x_3^{2^{2n}}&\dots&x_{2n}^{2^{2n}}&0\\
\end{matrix} 
\right)
\left(
\begin{matrix} x_2&x_2^2&x_2^4&\dots&x_2^{2^{2n}}\\\\
x_1&x_1^2&x_1^4&\dots&x_1^{2^{2n}}\\\\
x_4&x_4^2&x_4^8&\dots&x_4^{2^{2n}}\\\\
\vdots&\vdots&\vdots&\ddots&\vdots\\\\
x_{2n-1}&x_{2n-1}^2&x_{2n-1}^4&\dots&x_{2n-1}^{2^{2n}}\\\\
0&0&0&\dots&0\\
\end{matrix} \right)$$
$$=\alpha\left( \begin{matrix} 2X & \xi_1 & \xi_2 & \xi_3 & \dots & \xi_{2n}
\\\\ 
\xi_1 & 2X^2 & \xi_1^2 & \xi_2^2 & \dots & \xi_{2n-1}^2 \\\\ 
\xi_2 &
\xi_1^2 & 2X^4 & \xi_1^4 & \dots & \xi_{2n-2}^4 \\\\ 
\xi_3 & \xi_2^2 &
\xi_1^4 & 2X^8 & \dots & \xi_{2n-3}^8 \\\\ 
\vdots &\vdots &\vdots &\vdots&\ddots&\vdots\\\\ 
\xi_{2n} & \xi_{2n-1} ^2 & \xi_{2n-2} ^4 & \xi_{2n-3} ^8
& \dots & 2X^{2^{2n}} \\ \end{matrix} \right).
$$
On taking determinants we find that $\alpha\left(H_{2n}\right)=0$ and hence $\alpha\left(\frac12H_{2n}\right)=0$. (Recall from the remarks preceding Definition \ref{om def} that $H_{2n}$ is divisible by $2$.)
By definition, $\Omega_{2n}(X) := (\frac{1}{2} \; H_{2n})$ mod $2$. Now the image of 
$\alpha\left(\frac12H_{2n}\right)$ under the map $\Z\to\F_2$ is $\Omega(\xi_0+x_0^2)$ and hence $\Omega(\xi_0+x_0^2)=0$ as required.
\item
By part (ii), $\Omega_{2n}=c_0^2Q(X)=c_0^2Q^-(X)Q^+(X)$. 
The 
Dickson element $c_0$ and the quadratic Chern polynomials $Q^\pm(X)$ are irreducible as the symplectic group transitively permutes their factors.
\item
Lemma \ref{deg arg}(ii) deals with the case of $c_0Q^-(X)$. Part (iv) of that Lemma says that $c_0^2Q^+(X)$ belongs to our target ring.
Thus we have 
\begin{eqnarray*}
c_0^3Q^+(X)Q^-(X) &=& c_0^2Q^+(X) \cdot c_0Q^-(X) \\
&=& c_0\cdot \Omega_{2n}(X),
\end{eqnarray*}
and this 2-way factorization happens in the polynomial ring $\F_2[X,\xi_1,\dots,\xi_{2n}]$. In this ring, $c_0$ is prime. Hence either 
$c_0$ divides $c_0Q^-(X)$ or $c_0$ divides $c_0^2Q^+(X)$ in $\F_2[X,\xi_1,\dots,\xi_{2n}]$. 

Suppose that $c_0\,\vert \,c_0Q^-(X)$, so that $Q^-(X)\in \F_2[X,\xi_1,\dots,\xi_{2n}]$. By Lemma \ref{p squared} and Lemma \ref{mis1}, the restriction of $Q^-(X)$ to a maximal isotropic subspace $W$ of $V$ involves Dickson invariants for $W$. As each $\xi_i$ restricts to zero on $W$, it follows that $Q^-(X)$ involves Dickson invariants of $U$. This contradicts the assumption, so it must be the case that 
$c_0\,\vert \,c_0^2Q^+(X)$ in $\F_2[X,\xi_1,\dots,\xi_{2n}]$, and $c_0Q^+(X)\in \F_2[X,\xi_1,\dots,\xi_{2n}]$
 as claimed. 
\end{enumerate} 
\end{proof}

\begin{definition}\label{defompm}
Motivated by the previous Lemma, we define polynomials 
$\Omega^+_{2n}(X)$ and $\Omega^+_{2n}(X)$ in $\F_2[X,\xi_1,\dots,\xi_{2n}]$ by
$$
\Omega^+_{2n}(X):=c_0 Q^+(X), \quad \Omega^-_{2n}(X):=c_0 Q^-(X).
$$
\end{definition}
Part (iii) of the Lemma says that $\Omega_{2n}(X)=\Omega^+_{2n}(X)\Omega^-_{2n}(X)$.

\begin{lemma}These polynomials have the following properties.
\begin{enumerate} 
\item 
$\Omega_{2n}^-(X)$ and
$\Omega_{2n}^+(X)$ are irreducible in the polynomial ring $\F_2[X,\xi_1,\dots,\xi_{2n}]$.
\item 
$\Omega^\pm_{2n}(X)$ are both linear in $\xi_{2n}$ with coefficients 
$\left(\Omega_{2n-2}^{\pm}(X)\right)^2$. 
\end{enumerate}
\end{lemma}

\begin{proof}
\begin{enumerate}
\item
Lemma \ref{omegas} part (iv) says that $\Omega_{2n}^{\pm}(X)\in\F_2[X,\xi_1,\dots,\xi_{2n}]$, 
and part (iii) that   
 $c_0$, $Q^+(X)$ and $Q^-(X)$ are irreducible elements of the ring 
 $$S(U^*)^{Sp(U)}[X]=\F_2[X,\xi_1,\dots,\xi_{2n},c_{2n-1},\dots, c_{n}].$$
Thus if $\Omega_{2n}^{\pm}(X)$ were reducible in the smaller ring $\F_2[X,\xi_1,\dots,\xi_{2n}]$ it would factorize as $c_0\cdot Q^{\pm}(X)$. By the argument in Lemma \ref{omegas} part (iv) we know however that $Q^{\pm}(X)$ is not in $\F_2[X,\xi_1,\dots,\xi_{2n}]$, so $\Omega_{2n}^{\pm}(X)$ is irreducible in that ring as claimed. 

\item View $\Omega_{2n}(X)$ as a polynomial in $\xi_{2n}$ with coefficients in $\F_2[X,\xi_1,\dots,\xi_{2n-1}]$.
From the Definition \ref{om def} of $\Omega_{2n}(X)$ via the determinant $H_{2n}$ we see that $\Omega_{2n}(X)$ is quadratic in $\xi_{2n}$ and
$$
\Omega_{2n}(X) = \xi_{2n}^2\Omega_{2n-2}(X)^2 + \xi_{2n}\text{-linear terms.}$$
By Lemma \ref{omegas} part (iii), the quadratic term is $\xi_{2n}^2(\Omega_{2n-2}^-(X))^2 (\Omega_{2n-2}^+(X))^2$.
By part (i), each of   $\xi_{2n}$, $\Omega_{2n-2}^+(X)$ and $\Omega_{2n-2}^-(X)$ is irreducible in the polynomial ring $\F_2[X,\xi_1,\dots,\xi_{2n}]$. Further, 
$\Omega_{2n}(X)=\Omega^+_{2n}(X)\Omega^-_{2n}(X)$, and a degree argument now delivers the result.
\end{enumerate}
\end{proof}

\begin{lemma}\label{omegaidentity}
\begin{enumerate} 
\item If $q$ is any quadratic form of
$-$type then 
\begin{eqnarray*} \Omega_{2n}^-(q)&=&0\\
\Omega_{2n}^+(q)&=&c_0\prod_{q+x^2 \text{ has $+$type}}x^2\\
\Omega_{2n-2}^-(q)&=&\prod_{x\ne0\ \&\ q+x^2 \text{ has $-$type}}x
\end{eqnarray*} 
\item If $q$ is any quadratic form of $+$type then
\begin{eqnarray*} 
\Omega_{2n}^-(q)&=&c_0\prod_{q+x^2 \text{ has
$-$type}}x^2\\ 
\Omega_{2n}^+(q)&=&0\\ 
\Omega_{2n-2}^+(q)&=&\prod_{x\ne0\
\&\ q+x^2 \text{ has $+$type}}x \end{eqnarray*} 
\item The following is
an identity:
$$\left(\Omega_{2n-2}^+(X)\right)^2\Omega_{2n}^-(X)+\left(\Omega_{2n-2}^-
(X)\right)^2\Omega_{2n}^+(X)=\Lambda_{2n}^3.$$ \end{enumerate}
\end{lemma} 

\begin{proof}
\begin{enumerate} \item The first equality is clear. For the second, we
have 
$$
\Omega_{2n}^+(q):=c_0 \prod_{r\in B^+}(q+r).
$$
Any two of the singular quadratic forms differ by a square
in $U^*$, so each factor
$q+r$ equals $x^2$ for some $x\in U^*$. Then the product over the $r$ is the product over
the $q+x^2$ of $+$type, and the equality holds.

For the third equality, recalling our discussion in \S\ref{quadforms}, for any $q$ 
and any $x\neq 0$ the following are equivalent: 
\begin{itemize}
\item $q+x^2$ has $-$type
\item $q\vert_{\ker(x)}$ is a singular form of $-$type, 
\end{itemize}
and by the first equality this is equivalent to 
\begin{itemize}
\item $\Omega_{2n-2}^-(q)\vert_{\ker(x)}=0$. 
\end{itemize}

Consider the quotient space to $V$ obtained by restricting to the kernel
$\ker(x)$ of a non-zero element of $V^*$ of $-$type, i.e. an element for which
$q+x^2$ is of $-$type. $\Omega_{2n-2}^-(q)$ then restricts to zero on this
space by the remarks above. Thus $\Omega_{2n-2}^-(q)$ must be
divisible by $x\,$ as a polynomial. 

But this is true for any such $x$, so $\Omega_{2n-2}^-(q)$ is divisible
by their product. Comparing degrees now shows that it is in
fact equal to this product, as required.

\item Similar to part (i).

\item The left side of the equation is a polynomial with degree in $X$
at most equal to 
\begin{eqnarray*}
2\deg(\Omega_{2n-2}^+(X)) + \deg(\Omega_{2n}^-(X)) &=&  2(2^{2n-3} - 2^{n-2}) + (2^{2n-1} + 2^{n-1}) \\
&=& (2^{2n-1} + 2^{2n-2}) \\
 ( &=& 2\deg(\Omega_{2n-2}^-(X)) + \deg(\Omega_{2n}^+(X))\, )
\end{eqnarray*}
which is less than $2^{2n}$, the number of quadratic forms of $+$ or $-$
type. Thus if the equality holds when evaluated on 
every such form, it is an identity.

Now let $q$ be a form, say of $+$type. By parts (i) and (ii), the second
summand of the left hand side is zero, and the first is equal to
\begin{eqnarray*}
\left(\Omega_{2n-2}^-(q)\right)^2\Omega_{2n}^+(q) &=&
\left(\prod_{x\ne 0 \,\, \& \,\, q+x^2 \text{ has $+$type}}x \right)^2 
\left(c_0\prod_{q+x^2 \text{ has $-$type}}x^2 \right) \\
&=&
c_0\prod_{x\in U^*, \,\, x\ne 0}x^2 \\
&=& c_0^3 
\end{eqnarray*}
which is equal to $\Lambda_{2n}^3$ by Lemma \ref{lambdas} (i). 

The calculation is similar if $q$ is a form of $-$type, so as noted
above we are done.
\end{enumerate}
\end{proof}

This has important consequences for the Chern polynomials $P^\pm(t)$.

\begin{corollary}\label{coeff t in theps}
The coefficient of $t$ in $P^\pm(t)$ is $\Omega^\pm_{2n-2}(\xi_0)$.
\end{corollary}
\begin{proof}
We deal with the $-$case, the proof in the $+$case being obtained simply by replacing $+$ by $-$ throughout the argument.
Let $s$ be the coefficient of $t$ in $P^-(t)$. 
Let $x_-\in A^-$ be any vector of $-$type (i.e. $\xi_0+x_-^2$ has $-$type). On restricting to the subspace $\ker x_-$ we have
$$P^-(t)|_{\ker x_-}=t\cdot\prod_{x'}(t+x')$$ where $x'$ runs through the non-zero vectors in $\ker x_-$ such that $\xi_0+x_-^2+x'^2$ has $-$type. Equating coefficients of $t$ we see that
$$s|_{\ker x_-}=\prod_{x'}x'=\Omega^-_{2n-2}(\xi_0+x_-^2),$$ by Lemma \ref{omegaidentity}(i). Hence $s\equiv\Omega^-_{2n-2}(\xi_0)$ modulo $x_-$. This is true for all $x_-\in A^-$ and hence $s\equiv\Omega^-_{2n-2}(\xi_0)$ modulo $\prod x_-=P^-(0)$. Now $P^-(0)$ has degree greater than $s$, and so in fact we have
$$s=\Omega^-_{2n-2}(\xi_0)$$ as claimed.
\end{proof}

\begin{lemma} \label{11.6}
Viewing $\Omega_{2n}^\pm(X)$ as a polynomial in $X$:
\begin{enumerate} 
\item $\Omega_{2n}(X)=\left(\Lambda_{2n}\right)^2X^{2^{2n}}+\text{ terms involving lower powers of }X$. 
\item $\Omega_{2n}^+(X)=\Lambda_{2n}X^{2^{2n-1}+2^{n-1}}+\text{ terms involving lower powers of }X$. 
\item $\Omega_{2n}^-(X)=\Lambda_{2n}X^{2^{2n-1}-2^{n-1}}+\text{ terms involving lower powers of }X$. 
\end{enumerate}
\end{lemma} 

The proof of this is trivially a consequence of the Definition \ref{defompm}.

\begin{lemma} \label{nice formulae}
\begin{enumerate}
\item
$\Omega_{2n-2}^+(t^2+\xi_0)P^-(t)+\Omega_{2n-2}^-(t^2+\xi_0)P^+(t)=
\Lambda_{2n},$
\item
$\Lambda_{2n}P^-(t)=\Omega_{2n-2}^+(t^2+\xi_0)Q^-(t^2+\xi_0)+\Omega_{2n-2}^-
(t^2+\xi_0)\left(D(t)+D(x_0)\right),$ 
\item
$\Lambda_{2n}P^+(t)=\Omega_{2n-2}^-(t^2+\xi_0)Q^+(t^2+\xi_0)+\Omega_{2n-2}^+
(t^2+\xi_0)\left(D(t)+D(x_0)\right).$
\end{enumerate}
\end{lemma}
\begin{proof}
The first equation follows immediately from 
Lemma \ref{omegaidentity} (iii), 
putting $X=t^2 + \xi_0$, dividing by $\Lambda_{2n}$ and taking the
square root. The second follows from the first through multiplying by 
$P^-(t)$ and using the third equality of Lemma \ref{p squared}.
\end{proof}

We conclude this section with statements of our calculations of the $\Omega_{2n}^\pm$ for small values of $n$.

\begin{example}\label{eg omega}
For $n=1$ we have \begin{eqnarray*}
\Omega_2(X)&=&\xi_1^2X^4+\xi_2^2X^2+\xi_1^4X+\xi_1^3\xi_2\\
\Omega^+_2(X)&=&\xi_1X^3+\xi_2X^2+\xi_1^3\\ 
\Omega^-_2(X)&=&\xi_1X+\xi_2
\end{eqnarray*} 

For $n=2$,
\begin{eqnarray*}
\Omega_4(X)&=&\Lambda_4^2X^{16}+\Lambda_{4,3}^2X^8+\Lambda_{4,2}^2X^4
+\Lambda_{4,1}^2X^2+\Lambda_4^4X+\\
&&\quad (\xi_1^4\Lambda_{4,3}+(\xi_2^4+\xi_2\xi_1^5)\Lambda_4)
(\Lambda_{4,2}+\xi_1\xi_3\Lambda_4)\\
\Omega_4^+(X)&=&\Lambda_4X^{10}+\xi_1^2\Lambda_4X^7+(\Lambda_{4,3}+
\xi_1\xi_2\Lambda_4)X^6+\xi_2^2\Lambda_4X^5+\\
&&\quad(\Lambda_{4,2}+(\xi_1\xi_3+\xi_1^4)\Lambda_4)X^4+
\xi_2^2\xi_1^2\Lambda_4X^2+\xi_1^6\Lambda_4X+\\
&&\qquad\xi_1^4\Lambda_{4,3}+(\xi_2^4+\xi_2\xi_1^5)\Lambda_4\\
\Omega_4^-(X)&=&\Lambda_4X^6+\xi_1^2\Lambda_4X^3+(\Lambda_{4,3}+
\xi_1\xi_2\Lambda_4)X^2+\xi_2^2\Lambda_4X+(\Lambda_{4,2}+\xi_1\xi_3\Lambda_4)
\end{eqnarray*}

For $n=0$ it makes sense to define $\Omega^+_0(X):=X$ and $\Omega^-_0(X):=1$.
\end{example}

When $n=3$ the calculation is formidable. The factors $\Omega_6^-(X)$ and $\Omega_6^+(X)$ have degrees $119$ and $135$ respectively. The calculation was carried out directly from our definition of $\Omega_6(X)$ by Allan Steel using the Magma computer algebra package. His results are available on the second author's web site, \cite{steel}.
%%%%%%%%%%%%%%%%%%%%%%%%%%%%%%%%%%%%%%%%%%%%%%%%%%%%%%%%

\section{The recursive calculation of $\Omega^-(X)$ and $\Omega^+(X)$}

\begin{proposition}\label{alpha+}
In the abstract polynomial ring $\F_2[X,\xi_1,\xi_2,\dots]$ there is a sequence of polynomials $\alpha^+_n(X)$ for $n\ge0$ with the following properties:
\begin{enumerate}
\item For each $n\ge0$, $\alpha^+_n(X)$ belongs to 
$\F_2[X,\xi_1,\dots,\xi_{2n-1}]$ and
$$\Omega^+_{2n}(X)=\sum_{\ell=0}^n\Lambda_{2n,n+\ell}\left(\alpha^+_\ell(X)\right)^{2^{n-\ell}}.$$
\item $\alpha_n^+(X)$ is monic in $X$ with $X$-degree $2^{2n-1}+2^{n-1}$.
\item 
\begin{eqnarray*}
\alpha^+_0(X)&=&X\\
\alpha^+_1(X)&=&X^3+\xi_1^2\\
\alpha^+_2(X)&=&X^{10}+\xi_1^2X^7+\xi_2\xi_1X^6+\xi_2^2X^5+\\
&&\qquad(\xi_3\xi_1+\xi_1^4)X^4+\xi_2^2\xi_1^2X^2+\xi_1^6X+\xi_2^4+\xi_2\xi_1^5.
\end{eqnarray*}
\item $\alpha^+_n(X)\equiv\xi_n^{2^n}$ modulo $X,\xi_1,\dots,\xi_{n-1}$.
\item For $n\ge1$, the summand of $\alpha^+_n(X)$ comprising all those terms of odd degree in $X$ is equal to $\left(\Omega^+_{2n-2}(X)\right)^2X$.
\end{enumerate}
\end{proposition}
\begin{proof}
\begin{enumerate}
\item We construct the $\alpha^+_n(X)$ inductively. Notice first that $\Lambda_{2,1}=\xi_2$, 
$\Lambda_{2,2}=\Lambda_2=\xi_1$, and therefore
$$\Omega^+_2(X)=\xi_1X^3+\xi_2X^2+\xi_1^3=\Lambda_{2,1}(\alpha^+_0(X))^2+\Lambda_2\alpha^+_1(X)$$ confirming our formula when $n=1$. We could be a little more economical by observing that the formula is also consistent with the case $n=0$ given the definitions $\Omega^+_0(X)=X$ and $\Lambda_0=1$. Now suppose that $n\ge2$ and that the $\alpha^+_j(X)$ have been chosen for $j<n$ so that
$$\Omega^+_{2n-2}(X)=\sum_{\ell=0}^{n-1}\Lambda_{2n,n+\ell}\left(\alpha^+_\ell(X)\right)^{2^{n-\ell}}.$$ Then the coefficient of $\xi_{2n}$ in 
$f(X):=\Omega^+_{2n}(X)-\displaystyle\sum_{\ell=0}^{n-1}\Lambda_{2n,n+\ell}\left(\alpha^+_\ell(X)\right)^{2^{n-\ell}}$ is
$$\left(\Omega^+_{2n-2}(X)\right)^2-\sum_{\ell=0}^{n-1}\left(\Lambda_{2n-2,n+\ell-1}\right)^2\left(\alpha^+_\ell(X)\right)^{2^{n-\ell}}$$ which is zero by induction. Therefore $f(X)$ is a polynomial in $X,\xi_1,\dots,\xi_{2n-1}$ and since the coefficients of $f(X)$ as a polynomial in $X$ are divisible by $\Lambda_{2n}$ in the ambient ring $S$ it follows that $f(X)$ is also divisible by $\Lambda_{2n}$ in the ring
$\F_2[X,\xi_1,\dots,\xi_{2n-1}]$. Therefore we can (and must) set $\alpha^+_n(X):=f(X)/\Lambda_{2n}$.
\item This follows from Lemma \ref{11.6}. Working by induction on $n$ we may suppose that $\alpha^+_\ell(X)$ is monic of $X$-degree $2^{2\ell-1}+2^{\ell-1}$ for $\ell<n$. Lemma \ref{11.6} says that $\Omega^+_{2n}(X)$ has $X$-degree $2^{2n-1}+2^{n-1}$ which is greater than any of the contributions coming from $\left(\alpha^+_\ell(X)\right)^{2^{n-\ell}}$ for $\ell<n$ in our new formula. Therefore only the term involving $\alpha^+_n(X)$ contributes to the highest power of $X$ and Lemma \ref{11.6} verifies our assertion that $\alpha^+_n(X)$ is monic.
\item This can be checked by direct calculation.
\item  The total degree of $\alpha_n^+(X)$ of is $2^{2n}+2^n$ because $X$ has degree $2$. Modulo $X,\xi_1,\dots,\xi_{n-1}$, only monomials in $\xi_n,\dots,\xi_{2n-1}$ can contribute. For such a monomial
$\xi_n^{a_n}\cdots\xi_{2n-1}^{a_{2n-1}}$ to exist we need natural numbers $a_n,\dots,a_{2n-1}$ such that 
$$a_n(2^n+1)+\dots+a_{2n-1}(2^{2n-1}+1)=2^{2n}+2^n.$$
Hence $$a_n+\dots+a_{2n-1}\equiv0\text{ modulo }2^n.$$ Moreover we can also see that $a_j\le2^{n-j}$ for each $j$ and therefore $$a_n+\dots+a_{2n-1}\le2^n+2^{n-1}+\dots+2<2^{n+1}.$$ If the $a_i$ are not all zero then the displayed information above implies that $$a_n+\dots+a_{2n-1}=2^n,$$
and therefore 
\begin{eqnarray*}
2^n(2^n+1)&=&(a_n+\dots+a_{2n-1})(2^n+1)\\
&\le&a_n(2^n+1)+\dots+a_{2n-1}(2^{2n-1}+1)\\
&=&2^{2n}+2^n.
\end{eqnarray*}
Thus the inequality above must be equality. Hence $a_i=0$ for $i>n$ and it follows that $a_n=2^n$ and $\xi_n^{a_n}\cdots\xi_{2n-1}^{a_{2n-1}}=\xi_n^{2^n}$.
From this we can conclude that $\alpha^+_n(X)\equiv0$ or $\xi_n^{2^n}$ modulo $X,\xi_1,\dots,\xi_{n-1}$. In particular, it follows that $\alpha^+_\ell(X)\equiv0$ modulo $X,\xi_1,\dots,\xi_{n-1}$ whenever $\ell<n-1$, and so the recursive formula in (i) gives
$$\Omega^+_{2n}(0)\equiv\Lambda_{2n}\alpha^+_n(X)\mod X,\xi_1,\dots,\xi_{n-1}.$$
It can be seen directly from its definition that $\Omega_{2n}(X)$ involves the monomial $\xi_{n+1}^{2^n-1}\xi_n^{2^{n+1}-1}$ and therefore $\Omega_{2n}(X)\not\equiv0$, $\Omega^+_{2n}(X)\not\equiv0$, and $\alpha^+_n(X)\not\equiv0$. Thus 
$\alpha^+_n(X)\equiv\xi_n^{2^n}$ modulo $X,\xi_1,\dots,\xi_{n-1}$ as claimed.
\item We leave the proof of this part as an exercise.
\end{enumerate}
\end{proof}

\begin{proposition}
In the abstract polynomial ring $\F_2[X,\xi_1,\xi_2,\dots]$ there is a sequence of polynomials 
$\alpha^-_n(X)$ for $n\ge0$ with the following properties:
\begin{enumerate}
\item For each $n\ge0$, $\alpha^-_n(X)$ belongs to 
$\F_2[X,\xi_1,\dots,\xi_{2n-1}]$ and
$$\Omega^-_{2n}(X)=\sum_{\ell=0}^n\Lambda_{2n,n+\ell}\left(\alpha^-_\ell(X)\right)^{2^{n-\ell}}.$$
\item $\alpha_n^-(X)$ is monic in $X$ with $X$-degree $2^{2n-1}-2^{n-1}$.
\item 
\begin{eqnarray*}
\alpha^-_0(X)&=&1\\
\alpha^-_1(X)&=&X\\
\alpha^-_2(X)&=&X^6+\xi_1^2X^3+\xi_2\xi_1X^2+\xi_2^2X+\xi_3\xi_1.
\end{eqnarray*}
\item For $n\ge1$, $\alpha^-_n(X)\equiv0$ modulo $X,\xi_1,\dots,\xi_{n-1}$.
\item For $n\ge1$, the part of $\alpha^-_n$ which has odd degree in $X$ is equal to $\left(\Omega^-_{2n-2}(X)\right)^2X$.
\end{enumerate}
\end{proposition}
\begin{proof}
The proof proceeds in just the same way as for the $\alpha^+$. We leave the details to the reader. Just one remark: the degree argument for part (ii) goes along rather more swiftly in this case but notice that while the conclusion is stronger for $n\ge1$ the case $n=0$ is a significant exception to the rule.
\end{proof}

\begin{corollary}\label{middle of omega}
For each $n\ge1$, we have that $\Omega^+_{2n}(X)\equiv\xi_n^{2^{n+1}-1}$ and
$\Omega^-_{2n}(X)\equiv\Lambda_{2n,n}$ modulo $X,\xi_1,\dots,\xi_{n-1}$. Also, 
$\Omega^-_{2n}(X)\equiv\xi_{n+1}^{2^n-1}$ modulo $\xi_1,\dots,\xi_n$.
\end{corollary}
\begin{proof}
Using the above Propositions we find that
$$\Omega^+_{2n}(X)\equiv\Lambda_{2n}\alpha^+_n(X)\equiv\Lambda_{2n}\xi_n^{2^n}$$
and
$$\Omega^-_{2n}(X)\equiv\Lambda_{2n,n}\alpha^-_0(X)=\Lambda_{2n,n}$$
modulo $X,\xi_1,\dots,\xi_{n-1}$. 
This demonstrates the result for $\Omega^-$. 
It is also straightforward to see from its definition that
$\Lambda_{2n}\equiv\xi_n^{2^n-1}$ and the result for $\Omega^+$ now follows as well.
Working modulo $\xi_1,\dots,\xi_n$ for the last part, it is straightforward to use degree arguments to see that $\Lambda_{2n,n+j}\equiv0$ for $j\ge1$ and by direct calculation that $\Lambda_{2n,n}\equiv\xi_{n+1}^{2^n-1}$. These facts yield the third stated equivalence.
\end{proof}

\begin{corollary}\label{caroline} Let $j\ge0$ be a natural number and let $P_j$ be the coefficient of degree $j$ in $P^-(t)$, (i.e. $P_j$ is the coefficient of $t^{2^{2n-1}-2^{n-1}-j}$). 
Let $P_j'$ be the coefficient of degree $j-1$ in $\Omega^-_{2n-2}(t^2+\xi_0)$, (i.e. $P_j'$ is the coefficient of $t^{2^{2n-1}-2^{n-1}-j}$).
Then
$$Sq^{j}\Omega^-_{2n-2}(\xi_0)=\Omega^-_{2n-2}(\xi_0)P_j+P_j'P^-(0).$$
\end{corollary}
\begin{proof}
If $j=0$ then $P_0=1$ and we define $P_0':=0$. Assume that $j\ge1$.
Using Lemma \ref{coeff t in theps} together with Lemma \ref{insight} we have
$$Sq^j\Omega^-_{2n-2}(\xi_0)\equiv\Omega^-_{2n-2}(\xi_0)P_j$$ modulo $x_-$ for all $x_-\in A^-$. Hence $P^-(0)=\prod x_-$ divides $Sq^j\Omega^-_{2n-2}(\xi_0)+\Omega^-_{2n-2}(\xi_0)P_j$ and there is an invariant $P_j'$ such that
$$Sq^{j}\Omega^-_{2n-2}(\xi_0)=\Omega^-_{2n-2}(\xi_0)P_j+P_j'P^-(0).$$
Since $\Omega^-_{2n-2}(\xi_0)$ involves only $\xi_0,\dots,\xi_{2n-2}$, when we apply the Steenrod operation $Sq^j$ the result is a polynomial in $\xi_0,\dots,\xi_{2n-1}$. Therefore if we square the above displayed equation and multiply by $\Lambda_{2n}$ the left hand side does not involve $\xi_{2n}$ whilst the right hand side 
simplifies to the expression
\begin{eqnarray*}
&&(\Omega^-_{2n-2}(\xi_0))^2\Lambda_{2n}(P_j)^2+(P_j')^2\Lambda_{2n}(P^-(0))^2\\
&=&(\Omega^-_{2n-2}(\xi_0))^2\left[\sum_{\ell=0}^n\Lambda_{2n,n+\ell}\left(\alpha^-_\ell(t^2+\xi_0)\right)^{2^{n-\ell}}\right]_{[2j+2^{2n}-1]}
+(P_j')^2\Omega^-_{2n}(\xi_0),
\end{eqnarray*}
where the notation $[\ \ ]_{[2j+2^{2n}-1]}$ means ``pick out the coefficient of degree $2j+2^{2n}-1 $'' from the polynomial in $t$. Thus, equating coefficients of $\xi_{2n}$ we have
\begin{eqnarray*}
0&=&(\Omega^-_{2n-2}(\xi_0))^2\left[\sum_{\ell=0}^{n-1}\left(\Lambda_{2n-2,n+\ell-1}\right)^2\left(\alpha^-_\ell(t^2+\xi_0)\right)^{2^{n-\ell}}\right]_{[2j-2]}
+(P_j')^2\left(\Omega^-_{2n-2}(\xi_0)\right)^2\\
&=&\left(\Omega^-_{2n-2}(\xi_0)\left[\Omega^-_{2n-2}(t^2+\xi_0)\right]_{[j-1]}+P_j'\Omega^-_{2n-2}(\xi_0)\right)^2
\end{eqnarray*}
Taking square roots and dividing by $\Omega^-_{2n-2}(\xi_0)$ we conclude that
$$P_j'=\left[\Omega^-_{2n-2}(t^2+\xi_0)\right]_{[j-1]}$$ as required.
\end{proof}
\begin{corollary}\label{carolinex}
$\left(\Omega^-_{2n-4}(\xi_0)\right)^2\xi_{2n-1}+\Omega^-_{2n-2}(\xi_0)d_{2n-1}+\Lambda_{2n-2}d_n$ belongs to the subring generated by $\xi_0,\dots,\xi_{2n-2}$
\end{corollary}
\begin{proof}
Take $j:=2^{2n-2}$ in Corollary \ref{caroline}. Since 
$$\Omega^-_{2n-2}(X)=\xi_{2n-2}\left(\Omega^-_{2n-4}(X)\right)^2+
\text{ terms in }X,\xi_1,\dots,\xi_{2n-3}$$ it follows from Lemma \ref{squarexi} and its Corollary that
$$Sq^{2^{2n-2}}\left(\Omega^-_{2n-2}(X)\right)
=\xi_{2n-1}\left(\Omega^-_{2n-4}(X)\right)^2+
\text{ terms in }X,\xi_1,\dots,\xi_{2n-2}.$$ Therefore by Corollary \ref{caroline} we have
$$\xi_{2n-1}\left(\Omega^-_{2n-4}(X)\right)^2+
\text{ terms in }X,\xi_1,\dots,\xi_{2n-2}=
\Omega^-_{2n-2}(\xi_0)P_{2^{2n-2}}+P_{2^{2n-2}}'d_n,$$ where $P_{2^{2n-2}}'$ is the coefficient of degree $2^{2n-2}-1$ in $\Omega^-_{2n-2}(t^2+\xi_0)$. Using Lemma \ref{11.6} we can see that $P_{2^{2n-2}}'=\Lambda_{2n-2}$ and the desired conclusion follows.
\end{proof}
\begin{corollary}
The coefficients of $\Omega^-_{2n-2}(\xi_0)P^-(t)$ lie in $\F_2[\xi_0,\dots,\xi_{2n-1},d_n]$. 
The $2n+1$ elements $\xi_0,\dots,\xi_{2n-1},d_n$ are algebraically independent in $S$ and the $2n+2$ elements $\xi_0,\dots,\xi_{2n-1},d_n,d_{2n-1}$ satisfy a single relation in $S$ which is linear in $\xi_{2n-1},d_n,d_{2n-1}$.
\end{corollary}
\begin{proof}
The first part is immediate from Corollary \ref{caroline}. Now we know that $S$ is integral over the ring generated by the coefficients of $P^-(t)$ at least for $n\ge2$. It follows that $S$ is algebraic over $\F_2[\xi_0,\dots,\xi_{2n-1},d_n]$ and so on grounds of Krull dimension, $\xi_0,\dots,\xi_{2n-1},d_n$ must be algebraically independent. The last part now follows at once. The relation is as described in Corollary \ref{carolinex}
\end{proof}
\begin{corollary}\label{J2} Let $R$ denote the subring $\F_2[\xi_0,\dots,\xi_{2n-2}]$.
\begin{enumerate}
\item $\{s\in R;\ sd_{2n-1}\in\xi_{2n-1}R+d_nR\}=\Omega^-_{2n-2}(\xi_0)R,$
\item $\{s\in R;\ sd_{n}\in\xi_{2n-1}R+d_{2n-1}R\}=\Lambda_{2n-2}R,$
\item $\{s\in R;\ s\xi_{2n-1}\in d_{2n-1}R+d_nR\}=\Omega^-_{2n-4}(\xi_0)^2R.$
\end{enumerate}
\end{corollary}
\begin{proof}
This follows in a manner similar to the proof of Lemma \ref{J} using the fact that $\Omega^-_{2n-2}(\xi_0),\Lambda_{2n-2},\Omega^-_{2n-4}(\xi_0)$ are distinct irreducible elements of $\F_2[\xi_0,\dots,\xi_{2n-2}]$.
\end{proof}
%%%%%%%%%%%%%%%%%%%%%%%%%%%%%%%%%%%%%%%%%%%%%%%%%
We have the following consequence which is crucial in our calculations:
\begin{lemma}\label{frcd}
There is a matrix $\mathbf{J}_n$ and column vector $\mathbf{F}_n$ so that
$$\left(\begin{matrix}c_n\\\vdots\\c_{2n-1}\end{matrix}\right)=
\mathbf{J}_n^{*2}\left(\begin{matrix}d_n^2\\\vdots\\d_{2n-1}^2\end{matrix}\right)+\mathbf{F}_n.$$
Here, $\mathbf{J}_n^{*2}$ denotes the matrix whose entries are the squares of the entries of the matrix $\mathbf{J}_n$. The matrix $\mathbf{J}_n$ is upper uni-triangular with entries in the subring generated by $\xi_0,\dots,\xi_{2n-3}$ and $\mathbf{F}_n$ is a column of polynomials in $\xi_0,\dots,\xi_{2n-2}$.
\end{lemma}
\begin{proof}
From the formula 
$$P^-(t))^2=Q^-(t^2+\xi_0)=\sum_{\ell=0}^nc_{n+\ell}\left(\alpha^-_\ell(t^2+\xi_0)\right)^{2^{n-\ell}}$$ we see that on equating appropriate powers of $t$,
$$d_j^2=c_j+\text{ terms involving and linear in }c_{j+1},\dots,c_{2n-1},c_{2n}.$$
Note that $c_{2n}=1$. Moreover the coefficients of $c_k$ here are all squares for $j+1\le k\le 2n-1$. So there is an upper uni-triangular matrix $\mathbf{U}_n$ and a column $\mathbf{V}_n$ both having polynomial entries in the $\xi$'s such that
$$\left(\begin{matrix}d_n^2\\\vdots\\d_{2n-1}^2\end{matrix}\right)=
\mathbf{U}_n^{*2}\left(\begin{matrix}c_n\\\vdots\\c_{2n-1}\end{matrix}\right)+\mathbf{V}_n.$$ We obtain the required form of result by setting $\mathbf{J}_n:=\mathbf{U}_n^{-1}$ and then $\mathbf{F}_n:=\mathbf{J}_n^{*2}\mathbf{V}_n$.
\end{proof}
\begin{corollary}
The subring of $S$ generated by $\xi_0,\dots,\xi_{2n-1}$ together with $d_{2n-1},\dots,d_n$ contains all the coefficients of $P^-(t)$.
\end{corollary}
\begin{proof}
Let $P_j$ be the coefficient of $P^-(t)$ of degree $j$. The formula $$P^-(t)^2=\sum_{\ell=0}^nc_{n+\ell}\left(\alpha^-_\ell(t^2+\xi_0)\right)^{2^{n-\ell}}$$ shows that $P_j^2$ can be expressed as a linear combination of $c_n,\dots,c_{2n-1},c_{2n}=1$ with coefficients in the subring $\F_2[\xi_0,\dots,\xi_{2n-1}]$ and for $\ell<n$, the coefficient of $c_{n+\ell}$ in this expression is a square. From this, together with the Lemma above, we see that there is a linear combination 
$$\lambda_{2n-1}d_{2n-1}+\dots+\lambda_nd_n$$ with each $\lambda_i$ belonging to 
$\F_2[\xi_0,\dots,\xi_{2n-1}]$ such that
$$P_j^2=\left(\lambda_{2n-1}d_{2n-1}+\dots+\lambda_nd_n\right)^2+\lambda$$
for some $\lambda\in\F_2[\xi_0,\dots,\xi_{2n-1}]$. Clearly $\lambda$ is a square in the ambient ring $S$ and therefore it is a square in the subring $\F_2[\xi_0,\dots,\xi_{2n-1}]$
by Lemma \ref{squares}. So
$$P_j=\lambda_{2n-1}d_{2n-1}+\dots+\lambda_nd_n+\sqrt{\lambda},$$
as required.
\end{proof}
\begin{corollary}\label{ds gen +}
The coefficients of $P^+(t)$ also belong to the subring specified by the Corollary above.
\end{corollary}
\begin{proof}
This can be proved in just the same way, using the formula
$$P^+(t)^2=\sum_{\ell=0}^nc_{n+\ell}\left(\alpha^+_\ell(t^2+\xi_0)\right)^{2^{n-\ell}}.$$
\end{proof}

%%%%%%%%%%%%%%%%%%%%%%%%%%%%%%%%%%%%%%%%%%%%%%%%%%%%%%%%%%%%%%%%%%%%%%%
%%%%%%%%%%%%%%%%%%%%%%%%%%%%%%%%%%%%%%%%%%%%%%%%%%%%%%%%%%%%%%%%%%%%%%%

\section{The invariants of $O(V)$} 

Define $T$ to be the subring of $S$ generated by
$\xi_0,\xi_1,\dots,\xi_{2n-1},d_{2n-1},\dots,d_n$. This section is devoted to proving that $T$ is the ring of invariants of $O(V)$. If $n\ge2$ then the set $A^-$ of vectors in $V^*$ of $-$type spans $V^*$ and hence $S$ is an integral and separable extension of the subring generated by the coefficients of $P^-(t)$. When $n=0$ or $1$,  $A^-$ does not span $V^*$, but in these cases the $A^+$ spans $V^*$ and so at least we can say that $S$ is integral and separable over the subring generated by the coefficients of $P^+(t)$. The Corollaries following Lemma \ref{frcd} say that our chosen subring $T$ contains the coefficients of both $P^-(t)$ and $P^+(t)$. Therefore 

\begin{lemma}
$S$ is integral and separable over $T$ and the field of fractions $\fof(T)$ of $T$ is the fixed field $\fof(S)^{O(V)}$.
\end{lemma}
\begin{proof}
Only the remark about fields of fractions remains to be proved. 
We chose the generators of $T$ to be invariant and so $\fof(T)\subseteq \fof(S)^{O(V)}$. Let $G$ be the Galois group of the extension $\fof(T)\subseteq\fof(S)$. (The first part of this Lemma shows that the extension is Galois.) Then $G$ fixes $\xi_0\in\fof(T)$ and so $G\subseteq O(V)$. Also $\fof(T)\subseteq\fof(S)^{O(V)}$ and so $G\supseteq O(V)$. Therefore $G=O(V)$ and $\fof(T)=\fof(S)^{O(V)}$.
\end{proof}
From this Lemma we see that the integral closure of $T$ is the fixed ring for $O(V)$. In fact $T$ {\em is} integrally closed. To prove this we consider a presentation of $T$ as the quotient of a certain abstract polynomial ring $T^*$.

\begin{statement}{\bf Working in $T^*$.}\end{statement}

Let $T^*$ denote an abstract polynomial ring on $3n$ generators with weighted degrees in accordance with our chosen generators of $T$. 
We consider the surjection $T^*\to T$ defined by mapping the abstract generators to the corresponding generators of $T$.
We shall call the generators of $T^*$ by the obvious names:
$$\xi_0,\dots,\xi_{2n-1},d_{2n-1},\dots,d_n$$
in order to economize on notation. In practice this means keeping very clear the distinction between working in $T^*$ and working in $T$. We shall identify a regular sequence $r_1,\dots,r_{n-1}$ of elements of $T^*$ which lie in the kernel of the map $T^*\to T$. We shall therefore find that there is an induced map
$$T^*/(r_1,\dots,r_n)\to T.$$
The polynomials $\Omega^-_{2n-2}(\xi_0)$ and $\Lambda_{2n-2}$ are polynomials in the $\xi$'s which can be viewed as elements of $T^*$ in the obvious way. We shall show that their images in $T^*/(r_1,\dots,r_{n-1})$ satisfy the hypotheses of Proposition 1.1 and hence 
$T^*/(r_1,\dots,r_{n-1})$ is a unique factorization domain and the map 
$$T^*/(r_1,\dots,r_{n-1})\to T$$
is an isomorphism. Now clearly $T$ is contained in the ring of invariants $S^{O(V)}$ and since it is integrally closed one only has to check the elementary Galois theory to conclude that
$$T=S^{O(V)}.$$
Note that $T$ contains the Dickson invariants for $U^*$ and the coefficients of both $P^-(t)$ and $P^+(t)$. Therefore $T$ also contains $\eta=P^-(0)P^+(0)$ and so
$$S^{Sp(V)}\subset T.$$
This puts the Galois theory in place, and since any subgroup of $GL(V)$ which fixes the quadratic form $\xi_0$ is a subgroup of the orthogonal group with reach the desired conclusion.

We use the notation $\mathbf{M}^{*2}$ to indicate the matrix obtained from a matrix $\mathbf{M}$ by squaring all its entries. In case $\mathbf{M}=\mathbf{N}^{*2}$ for some $\mathbf{N}$ we also use the notation $\sqrt{\mathbf{M}}$ to denote the matrix $\mathbf{N}$ which is uniquely determined by $\mathbf{M}$ when it exists. We write $\mathbf{M}'$ for the matrix obtained from $\mathbf{M}$ by omitting the first row.

%:The subring of $S$ generated by $\xi_0,\dots,\xi_{2n}$ together with $d_{2n-1},\dots,d_n$.
\begin{statement}{\bf The subring of $S$ generated by $\xi_0,\dots,\xi_{2n}$ together with $d_{2n-1},\dots,d_n$.}\end{statement}

We begin by considering the subring of $S$ generated by $\xi_0,\dots,\xi_{2n},d_{2n-1},\dots,d_n$. This makes a total of $3n+1$ generators. Recall (\ref{frs}):
the fundamental relations for the symplectic invariants $S(U^*)^{Sp(U)}$
$$\left(\mathbf{L}_n\mathbf{K}_n+\mathbf{R}_n\right)
\left(\begin{matrix}c_n\\\vdots\\c_{2n-1}\end{matrix}\right)=
\left(\begin{matrix}\xi_{2n}\\\vdots\\\xi_{n+1}^{2^{n-1}}\end{matrix}\right)+\mathbf{L}_n\mathbf{E}_n$$
The following observations are significant:
\begin{lemma}\label{12.b}\ 
\begin{enumerate}
\item The $(n-1)\times(n-1)$ matrix obtained by omitting the first row and last column of $\mathbf{L}_n\mathbf{K}_n+\mathbf{R}_n$ is $\left(\mathbf{L}_{n-1}\mathbf{K}_{n-1}+\mathbf{R}_{n-1}\right)^{*2}$.
\item Working modulo $\xi_1,\dots,\xi_{n-1}$, we have $\mathbf{L}_n\equiv0$ and $\mathbf{R}_n$ is upper triangular with diagonal entries $\xi_n,\xi_n^2,\dots,\xi_n^{2^{n-1}}$.
\item $$\det\left(\mathbf{L}_n\mathbf{K}_n+\mathbf{R}_n\right)=\Lambda_{2n}$$
\end{enumerate}
\end{lemma}
\begin{proof}
\begin{enumerate}
\item This follows from the definitions and Lemma \ref{redundancy of cs}.
\item This is entirely straightforward.
\item Let $\delta$ denote the determinant. From (ii) we can deduce that $\delta\equiv\xi_n^{2^n-1}$ modulo $\xi_1,\dots,\xi_{n-1}$, and in particular it follows that $\delta$ is non-zero. Since our relations are homogeneous it follows that $\delta$ is a homogeneous polynomial, and we see that it has degree $2^{2n}-1$. Multiplying both sides of the matrix equation \ref{frs} by the matrix $\left(\mathbf{L}_n\mathbf{K}_n+\mathbf{R}_n\right)^{\text{cof}}$ of cofactors we see that
$$\left(\begin{matrix}\delta c_n\\\vdots\\\delta c_{2n-1}\end{matrix}\right)=
\left(\mathbf{L}_n\mathbf{K}_n+\mathbf{R}_n\right)^{\text{cof}}\left(
\left(\begin{matrix}\xi_{2n}\\\vdots\\\xi_{n+1}^{2^{n-1}}\end{matrix}\right)+\mathbf{L}_n\mathbf{E}_n\right)$$ This shows that $\delta$ belongs to the ideal
$J$ of Lemma \ref{J} and consequently, on grounds of degree, $\delta=\Lambda_{2n}$.
\end{enumerate}
\end{proof}

We also have Lemma \ref{frcd}, relating the Dickson invariants $c_{2n-1},\dots,c_n$ 
and the squares $d_{2n-1}^2,\dots,d_n^2$ of our fundamental orthogonal invariants. We use this in order to replace all the $c$'s with $d$'s. 
Using Lemma \ref{frcd}, we obtain the following matrix equation of relations:
$$\left(\mathbf{L}_n\mathbf{K}_n+
\mathbf{R}_n\right)\mathbf{J}_n^{*2}\left(\begin{matrix}d_n^2\\\vdots\\d_{2n-1}^2\end{matrix}\right)=
\left(\begin{matrix}\xi_{2n}\\\vdots\\\xi_{n+1}^{2^{n-1}}\end{matrix}\right)+\mathbf{L}_n\mathbf{E}_n+\left(\mathbf{L}_n\mathbf{K}_n+
\mathbf{R}_n\right)\mathbf{F}_n$$
At this stage, the first relation simply gives expression for $\xi_{2n}$ in terms of other generators. So we throw away this relation and throw away the redundant $\xi_{2n}$. 
In matrix form the situation can be summarized by omitting the first rows of chosen matrices. 
$$\left(\mathbf{L}'_n\mathbf{K}_n+
\mathbf{R}'_n\right)\mathbf{J}_n^{*2}\left(\begin{matrix}d_n^2\\\vdots\\d_{2n-1}^2\end{matrix}\right)=
\left(\begin{matrix}\xi_{2n-1}^2\\\vdots\\\xi_{n+1}^{2^{n-1}}\end{matrix}\right)+\mathbf{L}'_n\mathbf{E}_n+\left(\mathbf{L}'_n\mathbf{K}_n+
\mathbf{R}'_n\right)\mathbf{F}_n.$$
Notice that the only entries of $\mathbf{L}'_n$ which are not squares are the entries in the first column. Since $\mathbf{K}_n$ begins with a row of zeroes, the first column of $\mathbf{L}'_n$ makes no impact on the product $\mathbf{L}'_n\mathbf{K}_n$ and this matrix has square entries.
More precisely
we have
$$\mathbf{L}'_n\mathbf{K}_n=\left(\mathbf{L}_{n-1}^{*2}\mathbf{K}_{n-1}^{*2}\ \ \  \mathbf{L}_{n-1}^{*2}\mathbf{E}_{n-1}^{*2}\right)$$
 Moreover, $\mathbf{R}'_n$ also has square entries. Therefore the left hand side of our matrix equation consists entirely of squares. 
On the right hand side, the first vector comprises squares. We understand much less about the remaining vector
$$\mathbf{L}'_n\mathbf{E}_n+\left(\mathbf{L}'_n\mathbf{K}_n+
\mathbf{R}'_n\right)\mathbf{F}_n,$$
but it must of course consist of elements which are squares in the ambient ring $S$. Since this mysterious vector is a column of polynomials in $\xi_0,\dots,\xi_{2n-1}$ it follows that its entries are squares within the ring $\F_2[\xi_0,\dots,\xi_{2n-1}]$ by Lemma \ref{squares}(i). Hence we can take the square root of our matrix equation to obtain what we shall call
\begin{statement}\label{fro}{\bf The fundamental system of relations for $S^{O(V)}$:}\end{statement}
{\small $$\left(\left(\mathbf{L}_{n-1}\mathbf{K}_{n-1}\ \ \ \mathbf{L}_{n-1}\mathbf{E}_{n-1}\right)+\sqrt{\mathbf{R'_n}}\right)\mathbf{J}_n
\left(\begin{matrix}d_n\\\vdots\\d_{2n-1}\end{matrix}\right)=\left(\begin{matrix}\xi_{2n-1}\\\vdots\\\xi_{n+1}^{2^{n-2}}\end{matrix}\right)+
\sqrt{\mathbf{L}'_n\mathbf{E}_n+\left(\mathbf{L}'_n\mathbf{K}_n+
\mathbf{R}'_n\right)\mathbf{F}_n}$$}
We wish to reorganise this matrix equation in two different ways. First, we wish to describe the relations so that $d_{2n-1}$ appears on the right hand side. To this end, let $\mathbf{G}_{n-1}$ denote the last column of the matrix 
$$\left(\left(\mathbf{L}_{n-1}\mathbf{K}_{n-1}\ \ \ \mathbf{L}_{n-1}\mathbf{E}_{n-1}\right)+\sqrt{\mathbf{R'_n}}\right)\mathbf{J}_n$$ and let $\mathbf{S}_{n-1}$ denote the $(n-1)\times(n-1)$ matrix obtained by deleting this column. Then we can write
\begin{statement}\label{lhr}{\bf the left handed reorganisation:}\end{statement}
$$\mathbf{S}_{n-1}
\left(\begin{matrix}d_n\\\vdots\\d_{2n-2}\end{matrix}\right)=
\left(\begin{matrix}\xi_{2n-1}\\\vdots\\\xi_{n+1}^{2^{n-2}}\end{matrix}\right)+
\sqrt{\left(\mathbf{L}'_n\mathbf{E}_n+\left(\mathbf{L}'_n\mathbf{K}_n+
\mathbf{R}'_n\right)\mathbf{F}_n\right)}+\mathbf{G}_{n-1}d_{2n-1}$$ Since $\mathbf{J}_n$ is upper uni-triangular, and also the matrix obtained by deleting the last column of $\mathbf{R}'_n$ is $\mathbf{R}_{n-1}^{*2}$, it follows that 
$$\mathbf{S}_{n-1}=(\mathbf{L}_{n-1}\mathbf{K}_{n-1}+\mathbf{R}_{n-1})\mathbf{J}''_n,$$
where $\mathbf{J}''_n$ is the matrix obtained by deleting the last row and last column of $\mathbf{J}_n$.
Notice that $\mathbf{J}''_n$ is upper uni-triangular and therefore
\begin{lemma}
$$\det\mathbf{S}_{n-1}=\det\left(\mathbf{L}_{n-1}\mathbf{K}_{n-1}+\mathbf{R}_{n-1}\right)=\Lambda_{2n-2}$$
\end{lemma}
\noindent by Lemma \ref{12.b}(iii). Secondly we wish to describe the relations so that $d_n$ appears on the right hand side of the equation. To do this, let $\mathbf{H}_{n-1}$ denote the first column of the matrix
$$\left(\left(\mathbf{L}_{n-1}\mathbf{K}_{n-1}\ \ \ \mathbf{L}_{n-1}\mathbf{E}_{n-1}\right)+\sqrt{\mathbf{R'_n}}\right)\mathbf{J}_n$$
and let $\mathbf{T}_{n-1}$ denote the $(n-1)\times(n-1)$ matrix obtained by deleting this column. Then we can write
\begin{statement}\label{rhr}{\bf the right handed reorganization:}\end{statement}
$$\mathbf{T}_{n-1}
\left(\begin{matrix}d_{n+1}\\\vdots\\d_{2n-1}\end{matrix}\right)=
\left(\begin{matrix}\xi_{2n-1}\\\vdots\\\xi_{n+1}^{2^{n-2}}\end{matrix}\right)+
\sqrt{\left(\mathbf{L}'_n\mathbf{E}_n+\left(\mathbf{L}'_n\mathbf{K}_n+
\mathbf{R}'_n\right)\mathbf{F}_n\right)}+\mathbf{H}_{n-1}d_{n}.$$ We have the following result:
\begin{lemma}
$$\det\mathbf{T}_{n-1}=\Omega^-_{2n-2}(\xi_0).$$
\end{lemma}
\begin{proof}
Modulo $\xi_0,\dots,\xi_{n-1}$, $\mathbf{T}_{n-1}$ is upper triangular with diagonal entries $\xi_n,\xi_n^2,\dots,\xi_n^{2^{n-2}}$ and therefore $\det\mathbf{T}_{n-1}$ is non-zero. We also see that $\det\mathbf{T}_{n-1}$ is homogeneous of degree $(2^n+1)(2^{n-1}-1)=2^{2n-1}-2^{n-1}-1=\deg\Omega_{2n-2}^-(\xi_0)$. When we invert $\det\mathbf{T}_{n-1}$, we can solve the relations to give expressions for $d_{n+1},\dots,d_{2n-1}$ in terms of $\xi_0,\dots,\xi_{2n-1},d_n$. Therefore $\det\mathbf{T}_{n-1}$ belongs to the ideal of Corollary \ref{J2}(i) and hence the result follows.
\end{proof}
This concludes our investigation of the subring $T$.

\begin{statement}{\bf Working in $T^*$} \end{statement}

We consider an abstract polynomial ring in $3n$ generators:
$$T^*:=\F_2[\xi_0,\dots,\xi_{2n-1},d_{2n-1},\dots,d_n].$$ By using the same symbols $\xi_i$, $d_j$ as for the subring $T$ we risk great confusion. We shall  maintain a clear distinction between our work in  $T^*$ and $T$. The relations we found holding in $T$ can be interpreted as {\em relators} in $T^*$. There are $n-1$ of these relators and they are expressed by the column vector
{\footnotesize {$$\left(\begin{matrix}r_1\\\vdots\\r_{n-1}\end{matrix}\right)=
\left(\left(\mathbf{L}_{n-1}\mathbf{K}_{n-1}\ \ \ \mathbf{L}_{n-1}\mathbf{E}_{n-1}\right)+\sqrt{\mathbf{R'_n}}\right)\mathbf{J}_n
\left(\begin{matrix}d_n\\\vdots\\d_{2n-1}\end{matrix}\right)+
\left(\begin{matrix}\xi_{2n-1}\\\vdots\\\xi_{n+1}^{2^{n-2}}\end{matrix}\right)+
\sqrt{\mathbf{L}'_n\mathbf{E}_n+\left(\mathbf{L}'_n\mathbf{K}_n+
\mathbf{R}'_n\right)\mathbf{F}_n}.$$}}
Notice that since $\mathbf{L}'_n\equiv0$ modulo $\xi_1,\dots,\xi_{n-1}$ we therefore have the simplification
$$\left(\begin{matrix}r_1\\\vdots\\r_{n-1}\end{matrix}\right)\equiv
\left(\begin{matrix}\xi_{2n-1}\\\vdots\\\xi_{n+1}^{2^{n-2}}\end{matrix}\right)+
\sqrt{\mathbf{R}'_n\mathbf{F}_n},$$ modulo $d_n,\dots,d_{2n-1}, \xi_1,\dots,\xi_{n-1}$. Modulo $\xi_1,\dots,\xi_{n}$, the last row of $\mathbf{R}'_n$ is zero and hence
$$r_{n-1}\equiv\xi_{n+1}^{2^{n-2}}.$$ More generally for $j\ge1$, the last $j$ rows of $\mathbf{R}'_n$ are zero modulo $\xi_1,\dots,\xi_{n+j-1}$ and therefore
$$r_{n-j}\equiv\xi_{n+j}^{2^{n-j}}.$$
This has the crucial consequence that
\begin{lemma}\label{main regular sequence}
$$d_n,\dots,d_{2n-1},\xi_0, \xi_1,\dots,\xi_{n-1},\xi_n,r_{n-1},r_{n-2},\dots,r_{1}$$ is a regular sequence in $T^*$. 
\end{lemma}
\begin{proof}
We manage regular sequences of homogeneous elements in a graded commutative ring using three simple devices:
\begin{itemize}
\item Permuting the terms of a regular sequence yields a regular sequence. 
\item Replacing the last term of a regular sequence by a proper power yields a regular sequence. In view of the first device, we can in fact replace any term of a regular sequence by a proper power.
\item If $a_1,\dots,a_j,\dots,a_k$ is a regular sequence and $b$ is any ring element such that the ideals
$(a_1,\dots,a_{j-1},a_j)$ and $(a_1,\dots,a_{j-1},b)$ are equal then the sequence obtained by replacing $a_j$ by $b$ is a regular sequence.
\end{itemize}
These facts are all simple consequences of the definition that a sequence $a_1,\dots, a_k$ is regular if and only if each term is a non-zero-divisor modulo its predecessors. In $T^*$ we surely have the regular sequence
$$d_n,\dots,d_{2n-1},\xi_0, \xi_1,\dots,\xi_{n-1},\xi_n,\xi_{n+1},\xi_{n+2},\dots,\xi_{2n-1}$$
and since $r_1\equiv\xi_{2n-1}$ modulo $d_n,\dots,d_{2n-1},\xi_1,\dots,\xi_{2n-2}$ it follows that we can adjust the last term: so
$$d_n,\dots,d_{2n-1},\xi_0, \xi_1,\dots,\xi_{n-1},\xi_n,\xi_{n+1},\xi_{n+2},\dots,\xi_{2n-2},r_1$$
is a regular sequence. Therefore we can replace the penultimate term by its square and
$$d_n,\dots,d_{2n-1},\xi_0, \xi_1,\dots,\xi_{n-1},\xi_n,\xi_{n+1},\xi_{n+2},\dots,\xi_{2n-2}^2,r_1$$
is also a regular sequence. Since $r_2\equiv\xi_{2n-2}^2$ modulo $d_n,\dots,d_{2n-1},\xi_1,\dots,\xi_{2n-3}$ we see that
$$d_n,\dots,d_{2n-1},\xi_0, \xi_1,\dots,\xi_{n-1},\xi_n,\xi_{n+1},\xi_{n+2},\dots,\xi_{2n-3},r_2,r_1$$
is a regular sequence. Now we replace $\xi_{2n-3}$ by $\xi_{2n-3}^4$ and then by $r_3$. Continuing in this way the desired conclusion follows.
\end{proof}
The relators all belong to the kernel of the natural surjection $$T^*\to T.$$ Thus we have an induced map
$$T^*/(r_1,\dots,r_{n-1})\to T.$$
It is our aim to prove that this is an isomorphism and simultaneously that $T$ is integrally closed. 
We shall use Proposition 1.1. More precisely, we shall show that
\begin{lemma}
\begin{enumerate}
\item $\Lambda_{2n-2},\Omega^-_{2n-2}(\xi_0)$ is a regular sequence in 
$T^*/(r_1,\dots,r_{n-1})$;
\item the localizations 
$$T^*/(r_1,\dots,r_{n-1})[\Lambda_{2n-2}^{-1}],$$
$$T^*/(r_1,\dots,r_{n-1})[\Omega^-_{2n-2}(\xi_0)^{-1}]$$ are unique factorization domains;
\item $\Lambda_{2n-2}$ generates a prime ideal in  
$T^*/(r_1,\dots,r_{n-1})[\Omega^-_{2n-2}(\xi_0)^{-1}]$; and
\item $\Omega^-_{2n-2}(\xi_0)$ generates a prime ideal in  
$T^*/(r_1,\dots,r_{n-1})[\Lambda_{2n-2}^{-1}]$.
\end{enumerate}
\end{lemma}
\begin{proof} Note that to avoid excessive notation we are now in the position that the names of elements do not tell you which ring they belong to. When we write about the elements $\Lambda_{2n-2},\Omega^-_{2n-2}(\xi_0)$, keep in mind that the subring of $S$ generated by $\xi_0,\dots,\xi_{2n}$ is isomorphic to the abstract polynomial ring on $\xi_0,\dots,\xi_{2n}$ because these elements are algebraically independent in $S$. Here we are interested in viewing these polynomials in $T^*$ which is defined to be an abstract polynomial ring, and then in the quotient $T^*/(r_1,\dots,r_{n-1})$ which sits in between $T^*$ and $S$ as follows:
$$T^*\to T^*/(r_1,\dots,r_{n-1})\to T\subset S.$$
Therefore there is no real risk of confusion when considering a polynomial in the $\xi$'s so long as we keep clear which of the above rings is relevant at each point of argument.
\begin{enumerate}
\item We know from Lemma \ref{main regular sequence} that $$\xi_1,\dots,\xi_n$$ is a regular sequence in the quotient ring $T^*/(r_1,\dots,r_{n-1})$. Therefore
$$\xi_1,\dots,\xi_{n-1},\xi_n^{2^{n-1}-1}$$ is also a regular sequence. By Corollary \ref{middle of omega} we have that $\Omega^-_{2n-2}(\xi_0)\equiv\xi_{n}^{2^{n-1}-1}$ modulo $\xi_1,\dots,\xi_{n-1}$ and therefore
$$\xi_1,\dots,\xi_{n-1},\Omega^-_{2n-2}(\xi_0)$$ is a regular sequence, and so also is
$$\xi_1,\dots,\xi_{n-2},\xi_{n-1}^{2^{n-1}-1},\Omega^-_{2n-2}(\xi_0).$$ Using the fact that $\Lambda_{2n-2}\equiv\xi_{n-1}^{2^{n-1}-1}$ modulo $\xi_1,\dots,\xi_{n-2}$ we see that
$$\xi_0,\dots,\xi_{n-2},\Lambda_{2n-2},\Omega^-_{2n-2}(\xi_0)$$ is a regular sequence. In particular $$\Lambda_{2n-2},\Omega^-_{2n-2}(\xi_0)$$ is a regular sequence in $T^*/(r_1,\dots,r_{n-1})$ as claimed.
\item First notice from the left handed reorganization (\ref{lhr}) that when we invert $$\det\mathbf{S}_{n-1}=\Lambda_{2n-2}$$ the relations can be solved to express $d_n,\dots,d_{2n-2}$ in terms of $\xi_0,\dots,\xi_{2n-1}$ and 
$d_{2n-1}$. This gives the isomorphism 
$$T^*/(r_1,\dots,r_{n-1})[\Lambda_{2n-2}^{-1}]
\iso\F_2[\xi_0,\dots,\xi_{2n-1},d_{2n-1},\Lambda_{2n-2}^{-1}]$$ and we can see on grounds of Krull dimension that this ring is a localized polynomial ring; in particular it is a unique factorization domain.
Similarly, from the right handed reorganization (\ref{rhr}) we deduce the isomorphism
$$T^*/(r_1,\dots,r_{n-1})[\Omega^-_{2n-2}(\xi_0)^{-1}]
\iso\F_2[\xi_0,\dots,\xi_{2n-1},d_{n},\Omega^-_{2n-2}(\xi_0)^{-1}],$$
because inverting $\det\mathbf{T}_{n-1}=\Omega^-_{2n-2}(\xi_0)$ allows us to solve for $d_{n+1},\dots,d_{2n-1}$ in terms of $\xi_0,\dots,\xi_{2n-1},d_n$,
 and this ring is also a localized polynomial ring; in particular it is a unique factorization domain.
\end{enumerate}
Finally (iii) and (iv) hold because $\Lambda_{2n-2}^{-1}$ and $\Omega^-_{2n-2}(\xi_0)$ are irreducible polynomials in the ring $\F_2[\xi_0,\dots,\xi_{2n-1}]$.
\end{proof}
In view of Proposition 1.1, this completes the calculation of the ring of invariants $S^{O(V)}$.

%%%%%%%%%%%%%%%%%%%%%%%%%%%%%%%%%%%%%%%%%%%%%%%%%%%%%%%%%%%

\section{The invariants for $O^-$}

Let $\xi_-$ be any quadratic form of $-$type. Then $\xi_0+\xi_-=x_-^2$ for some fixed $x\in V^*\setminus U^*$. The composite $\ker x\to V\to U$ is an isomorphism and so we identify $\ker x$ with $U$. We write $O^-$ for the group of automorphisms of $(U,\xi_-)$

We show that $S(U^*)^{O^-}$ is given up to isomorphism by adding the single additional relation $P^-(0)=0$ to the ring of invariants $S^{O(V)}$. Let's see what happens when we adjoin this additional relation to our presentation. This simply amounts to setting $d_n=0$. 

Note that $\xi_-$ belongs to $S(U^*)$ and in just the same way as for Lemma \ref{6.3} we have
\begin{lemma}\label{minus0}
The elements $\xi_-,\xi_1,\dots,\xi_{2n-1}$ are algebraically independent in $S(U^*)$.
\end{lemma}
Moreover, Corollary \ref{caroline} simplifies in a significant way:
\begin{lemma}\label{caroline-}
For each $j$, $$Sq^j(\Omega^-_{2n-2}(\xi_-))=\Omega^-_{2n-2}(\xi_-)P_j$$
where $P_j$ is the coefficient of degree $j$ in the restriction of $P^-(t)$ to $U^*$.
\end{lemma}

We begin by studying the subring $T$ of $S(U^*)$ generated by the $3n-1$ elements $$\xi_-,\xi_1,\dots,\xi_{2n-1},d_{2n-1},\dots,d_{n+1}.$$ We shall in due course see that $T=S(U^*)^{O^-}$.

The fundamental system of relations for $S^{O(V)}$ now simplify in line with (\ref{rhr}):
\begin{statement}\label{minus1}
{\bf The fundamental system of relations for $O^-$.}
$$\mathbf{T}_{n-1}
\left(\begin{matrix}d_{n+1}\\\vdots\\d_{2n-1}\end{matrix}\right)=
\left(\begin{matrix}\xi_{2n-1}\\\vdots\\\xi_{n+1}^{2^{n-2}}\end{matrix}\right)+
\sqrt{\left(\mathbf{L}'_n\mathbf{E}_n+\left(\mathbf{L}'_n\mathbf{K}_n+
\mathbf{R}'_n\right)\mathbf{F}_n\right)}.$$ 
\end{statement}
Direct inspection of $\mathbf{T}_{n-1}$ reveals that its top right hand entry is of the form 
$$\xi_{2n-2}+\text{ terms involving }\xi_0,\dots,\xi_{2n-3}.$$ It follows that the $(n-2)\times(n-2)$ matrix $\mathbf{T}''_{n-1}$ found by deleting the first row and last column of $\mathbf{T}_{n-1}$ has determinant
$\Omega^-_{2n-4}(\xi_0)^2$. Since the first relation simply gives expression for $\xi_{2n-1}$ in terms of the other generators we can dispense with this relation, throw out $\xi_{2n-1}$, and use 
\begin{statement}\label{minus2}
{\bf the reduced system of relations for $O^-$.}
$$\mathbf{T}''_{n-1}
\left(\begin{matrix}d_{n+1}\\\vdots\\d_{2n-2}\end{matrix}\right)=
\left(\begin{matrix}\xi_{2n-2}\\\vdots\\\xi_{n+1}^{2^{n-2}}\end{matrix}\right)+
\sqrt{\left(\mathbf{L}'_n\mathbf{E}_n+\left(\mathbf{L}'_n\mathbf{K}_n+
\mathbf{R}'_n\right)\mathbf{F}_n\right)}+\mathbf{G}'_{n-1}d_{2n-1}.$$
\end{statement}

This concludes our analysis of the subring $T$. Now we consider an abstract polynomial ring on generators 
$$\xi_-,\xi_1,\dots,\xi_{2n-1},d_{2n-1},\dots,d_{n+1}.$$
The matrix relation (\ref{minus1}) can be interpreted as providing a sequence of relators $r_1,\dots,r_{n-1}$ in $T^*$ so that there is a natural map
$$T^*/(r_1,\dots,r_{n-1})\to T.$$
As before, $r_1,\dots,r_{n-1}$ form part of a longer regular sequence. This time 
\begin{statement}\label{minus3}
{\bf the fundamental regular sequence is}
$$d_{n+1},\dots,d_{2n-1},\xi_-,\xi_1,\dots,\xi_n,r_{n-1},\dots,r_1.$$
\end{statement}
In a simple variation on the odd dimensional case we have
\begin{lemma}
\begin{enumerate}
\item The sequence $$\Omega^-_{2n-4}(\xi_-),\Omega^-_{2n-2}(\xi_-)$$ is a regular sequence in $T^*/(r_1,\dots,r_{n-1})$,
\item the localizations
$$T^*/(r_1,\dots,r_{n-1})[\Omega^-_{2n-4}(\xi_-)^{-1}]$$
$$T^*/(r_1,\dots,r_{n-1})[\Omega^-_{2n-2}(\xi_-)^{-1}]$$
are isomorphic to the localized polynomial rings
$$\F_2[\xi_-,\xi_1,\dots,\xi_{2n-2},d_{2n-1},\Omega^-_{2n-4}(\xi_-)^{-1}]$$
$$\F_2[\xi_-,\xi_1,\dots,\xi_{2n-2},\xi_{2n-1},\Omega^-_{2n-2}(\xi_-)^{-1}]$$
respectively, and so these are unique factorization domains,
\item $\Omega^-_{2n-4}(\xi_-)$ is irreducible in 
$$\F_2[\xi_-,\xi_1,\dots,\xi_{2n-2},\xi_{2n-1},\Omega^-_{2n-2}(\xi_-)^{-1}]$$
and $\Omega^-_{2n-2}(\xi_-)$ is irreducible in 
$$\F_2[\xi_-,\xi_1,\dots,\xi_{2n-2},d_{2n-1},\Omega^-_{2n-4}(\xi_-)^{-1}].$$
\end{enumerate}
\end{lemma}
Proposition 1.1 now comes into force and we deduce that
$T^*/(r_1,\dots,r_{n-1})$ is isomorphic to the subring $T$ and that this is the ring of invariants as required.

%%%%%%%%%%%%%%%%%%%
%%%%%%%%%%%%%%%%%%%

\section{The invariants for $O^+$}

It is useful now to have formulations of Corollaries \ref{caroline} and \ref{carolinex} for $P^+(t)$. These are as follows and are proved in exactly the same way.

\begin{corollary} Let $j\ge0$ be a natural number and let $P^+_j$ be the coefficient of degree $j$ in $P^+(t)$, (i.e. $P^+_j$ is the coefficient of $t^{2^{2n-1}+2^{n-1}-j}$). 
Let $P_j''$ be the coefficient of degree $j-1$ in $\Omega^+_{2n-2}(t^2+\xi_0)$, (i.e. $P_j''$ is the coefficient of $t^{2^{2n-1}+2^{n-1}-j}$).
Then
$$Sq^{j}\Omega^+_{2n-2}(\xi_0)=\Omega^+_{2n-2}(\xi_0)P^+_j+P_j''P^+(0).$$
\end{corollary}

\begin{corollary}
$\left(\Omega^+_{2n-4}(\xi_0)\right)^2\xi_{2n-1}+\Omega^+_{2n-2}(\xi_0)d_{2n-1}+\Lambda_{2n-2}P^+(0)$ belongs to the subring generated by $\xi_0,\dots,\xi_{2n-2}$
\end{corollary}

%\begin{corollary}
%The coefficients of $\Omega^+_{2n-2}(\xi_0)P^-(t)$ lie in $\F_2[\xi_0,\dots,\xi_{2n-1},P^+(0)]$. 
%The $2n+1$ elements $\xi_0,\dots,\xi_{2n-1},P^+(0)$ are algebraically independent in $S$ and the $2n+2$ elements $\xi_0,\dots,\xi_{2n-1},P^+(0),d_{2n-1}$ satisfy a single relation in $S$ which is linear in $\xi_{2n-1},P^+(0),d_{2n-1}$.
%\end{corollary}

We shall have need of the following embellishment of our discussion of the odd dimensional case.
\begin{lemma}\label{+}
There exist polynomials $f_0,\dots,f_n\in\F_2[\xi_0,\dots,\xi_{2n-1}]$ such that
$$P^+(0)=\sum_{\ell=0}^{n-1}f_\ell d_{n+\ell}+f_n$$ and which satisfy the conditions $f_0=\xi_0^{2^{n-1}}$, and in general for $1\le j\le n$ we have $f_j\equiv\xi_j^{2^{n-1}}$ modulo $\xi_0,\dots,\xi_{j-1}$.
\end{lemma}
\begin{proof}
Since $P^+(0)^2=Q^+(\xi_0)$ we see from Proposition \ref{alpha+} that
$$P^+(0)^2=
\sum_{\ell=0}^nc_{n+\ell}\left(\alpha^+_\ell(\xi_0)\right)^{2^{n-\ell}},$$ and in addition that modulo $\xi_0,\xi_1,\dots,\xi_{j-1}$ 
$$P^+(0)^2\equiv\xi_j^{2^n}c_{n+j}+\text{ terms involving Dickson invariants of lower degree}.$$ From this we can deduce that
$$P^+(0)\equiv\xi_j^{2^{n-1}}d_{n+j}+\text{ terms involving $d$'s of lower degree}.$$
and the result follows.
\end{proof}

We now choose a quadratic form $\xi_+$ on $U$ of $+$type. As in the $-$type case, there is a vector $x_+\in V^*\setminus U^*$ such that $\xi_+=\xi_0+x_+^2$. The composite of inclusion and restriction supplies an isomorphism between $\ker x_+$ and $U^*$ so we identify these two. On restriction to $\ker x_+$ there is a new relation, namely $P^+(0)=0$. 

The identification $\ker x_+=U^*$ gives us a map
$$S\to S(U^*)$$ which carries the invariants of $O(V)$ into the ring of invariants of $O^+$. 
Therefore we wish to stick with the choice of $d_j$ as key generators for the new ring of invariants even though it may at first appear that working with coefficients chosen from $P^+(t)$ would be more natural. In fact it makes very little difference because
of Lemma \ref{ds gen +}.
From Lemma \ref{+} we see that this is linear in the $d$'s with coefficients $f_0,\dots,f_{n-1}$ and we adjoin the row of $f$'s as an additional last row to the matrix
$$\left(\left(\mathbf{L}_{n-1}\mathbf{K}_{n-1}\ \ \ \mathbf{L}_{n-1}\mathbf{E}_{n-1}\right)+\sqrt{\mathbf{R'_n}}\right)\mathbf{J}_n$$ at the left of the fundamental relations for $O(V)$ in (\ref{fro}). This yields an $n\times n$ matrix $\mathbf{M}_n$ and allows us to write the relations we have found for the group $O^+$ in the form
$$\mathbf{M}_n
\left(\begin{matrix}d_n\\\vdots\\d_{2n-1}\end{matrix}\right)=\left(\begin{matrix}\xi_{2n-1}\\\vdots\\\xi_{n+1}^{2^{n-2}}\\\\\xi_n^{2^{n-1}}\end{matrix}\right)+
\left(\begin{matrix}\sqrt{\mathbf{L}'_n\mathbf{E}_n+\left(\mathbf{L}'_n\mathbf{K}_n+
\mathbf{R}'_n\right)\mathbf{F}_n}\\\\\xi_n^{2^{n-1}}+f_n\end{matrix}\right).$$
A variation on the arguments we have used before leads to the conclusion that
$$\det\mathbf{M}_n=\Omega^+_{2n-2}(\xi_+)$$
and that
$$\det\mathbf{M}''_n=\Omega^+_{2n-4}(\xi_+)^2$$
where $\mathbf{M}''_n$ is the matrix obtained by omitting the first row and last column of $\mathbf{M}_n$. The same regular sequence arguments can be applied, and the ring of invariants for $O^+$ is established.

\end{document}